\documentclass[11pt]{amsart}
\usepackage{amssymb,amsmath,amsthm}
\usepackage{enumerate}
\usepackage{graphicx}
\usepackage{epsfig}
\usepackage{xcolor}
\usepackage[all]{xy}
\usepackage{url}
\usepackage{mathrsfs}
\usepackage{bbm}
\usepackage[normalem]{ulem}

\DeclareFontFamily{U}{mathx}{}
\DeclareFontShape{U}{mathx}{m}{n}{<-> mathx10}{}
\DeclareSymbolFont{mathx}{U}{mathx}{m}{n}
\DeclareMathAccent{\widehat}{0}{mathx}{"70}
\DeclareMathAccent{\widecheck}{0}{mathx}{"71}

\renewcommand{\Re}{\operatorname{Re}}
\renewcommand{\Im}{\operatorname{Im}}

\newcommand\restr[2]{{
  \left.\kern-\nulldelimiterspace 
  #1 
  \vphantom{\big|} 
  \right|_{#2} 
  }}

\DeclareMathOperator{\spa}{span}
\DeclareMathOperator{\cb}{cb}

\DeclareMathOperator{\CB}{\mathcal{CB}}

\DeclareMathOperator{\amconv}{amconv}
\DeclareMathOperator{\aconv}{aconv}
\DeclareMathOperator{\conj}{conj}

\DeclareMathOperator{\MAX}{MAX}
\DeclareMathOperator{\Lip}{Lip}
\DeclareMathOperator{\diag}{diag}
\DeclareMathOperator{\inte}{int}

\newcommand{\Lipo}{{\Lip_0}}
\newcommand{\Lipcb}{{\Lip_0^{\cb}}}
\newcommand{\Fcb}{\mathcal{F}^{\cb}}
\newcommand{\F}{\mathcal{F}}

\newcommand{\n}[1]{ \left\|#1\right\| }
\newcommand{\bn}[1]{ \big\|#1\big\| }

\newcommand{\tn}[1]{{\left\vert\kern-0.25ex\left\vert\kern-0.25ex\left\vert #1 
    \right\vert\kern-0.25ex\right\vert\kern-0.25ex\right\vert}}
\newcommand{\N}{{\mathbb{N}}}
\newcommand{\R}{{\mathbb{R}}}
\newcommand{\C}{{\mathbb{C}}}

\newcommand{\K}{{\mathbb{K}}}

\newcommand{\B}{{\mathcal{B}}}

\newcommand{\pair}[2]{{\langle #1, #2 \rangle}}
\newcommand{\mpair}[2]{{\langle\langle #1, #2 \rangle\rangle}}

\newtheorem{theorem}{Theorem}[section]
\newtheorem{lemma}[theorem]{Lemma}

\newtheorem{corollary}[theorem]{Corollary}
\newtheorem{proposition}[theorem]{Proposition}

\newtheorem{examples}[theorem]{Examples}

\newtheorem{remark}[theorem]{Remark}

\theoremstyle{definition}
\newtheorem{definition}[theorem]{Definition}

\numberwithin{equation}{section}

\def\ker{{\rm ker\, }}

\begin{document}

\title[Linearization of completely Lipschitz functions]{Linearization of completely Lipschitz functions}

\author[J.A. Ch\'avez-Dom\'inguez]{Javier Alejandro Ch\'avez-Dom\'inguez}
\address{David and Judi Proctor Department of Mathematics, University of Oklahoma, Norman, OK 73019-3103,
USA} \email{jachavezd@ou.edu}

\author[V. Dimant]{Ver\'onica Dimant}
\address{Departamento de Matem\'{a}tica y Ciencias, Universidad de San
Andr\'{e}s, Vito Dumas 284, (B1644BID) Victoria, Buenos Aires,
Argentina and CONICET} \email{vero@udesa.edu.ar}

\thanks{The first-named author was partially supported by NSF grants DMS-2247374 and DMS-2555053. 
The second-named author was partially supported by CONICET PIP 11220250100329CO.
}


\date{}

\begin{abstract}
    We develop a linearization theory for completely Lipschitz mappings between operator spaces. Motivated by the classical theory of Lipschitz-free spaces, we introduce a completely Lipschitz-free space and show that completely Lipschitz mappings admit a canonical linearization as completely bounded operators on this space. We provide several equivalent descriptions of this operator space, extending to the noncommutative setting the classical constructions of Arens--Eells, de Leeuw, and Dixmier--Ng. In particular, we obtain a molecular description, a quotient representation of an operator $\ell_1$-space, and a realization as the closed linear span of evaluation functionals. Our framework applies over both the real and complex scalar fields; in the real case, we first establish some operator space results needed to carry out the constructions. We also show that our notion of completely Lipschitz mapping encompasses the previously introduced extensions of Lipschitz mappings to operator spaces. These results provide a unified framework for studying Lipschitz phenomena in the setting of operator spaces and establish the corresponding free operator spaces as a natural tool for their linearization.
\end{abstract}

\maketitle

\section{Introduction}

The linearization of Lipschitz functions (defined on a metric space and fixing a base point) \cite{Arens-Eells, Weaver-Lipschitz-algebras} is a well-established and fruitful tool in the study of nonlinear maps between metric and Banach spaces. It allows one to represent a Lipschitz mapping as a bounded linear operator acting on a suitable Banach space, thereby bringing the  techniques of linear functional analysis to bear on nonlinear problems. This perspective has led to the development of the theory of Lipschitz-free spaces, which provides a natural linear framework for Lipschitz mappings and has become an important object of study in its own right: after the introduction of Lipschitz-free spaces to Banach space theory by Godefroy and Kalton \cite{Godefroy-Kalton}, these spaces have commanded significant attention and interest in the field for the past two decades. Given the typical interaction between Banach space and operator space theories, where the latter is considered as a quantized or noncommutative version of the former, it is natural to ask to what extent these linearization procedures can be developed in the noncommutative setting.
Before one can even talk about linearization, an appropriate notion of Lipschitz function is needed.
A natural first attempt was the \emph{completely Lipschitz functions} of \cite{Braga-Sinclair}, but this definition turned out to be essentially linear: it was proved shortly afterwards \cite{Braga-CD} that such functions must be translations of real-linear functions.
A different approach was taken in \cite{Braga-CD-Sinclair}, which succeeded both in covering nonlinear functions and in having a corresponding linearization procedure, but with the caveat of not using the full operator space structure (specifically, only finitely many matricial levels).
Recently, \cite{braga2024small} succeeded in achieving the balance of a definition of Lipschitz function which admits nonlinear examples and uses all of the matricial levels.

The present paper was born out of the desire to develop a corresponding linearization theory for the Lipschitz functions of \cite{braga2024small}.
After doing so, we realized that everything works in a more general conceptual setup. Thus, our main contribution is to introduce a notion
of completely Lipschitz mappings \emph{relative} to a matrix set and study the corresponding linearization in the operator space framework.
This definition encompasses the aforementioned ones from
\cite{Braga-CD, braga2024small, Braga-CD-Sinclair} and also the classical notion of a Lipschitz map on a metric space; moreover, our general linearization result recovers both the one from \cite{Braga-CD-Sinclair} and the classical one on metric spaces.

While our definition (see Definition \ref{defn-completely-Lipschitz}) might seem daunting at first, its gist is easy to explain in the classical setting. Given a metric space $M$, the Lipschitz condition for a function $f : M \to \R$ is defined via $|f(x)-f(y)| \le C d(x,y)$ \emph{for all} pairs $(x,y)\in M\times M$. One can define a weaker notion of Lipschitz function by requiring the condition only on a distinguished subset of pairs $D \subseteq M \times M$, and our approach is an operator version of this simple idea. 

The theory of Lipschitz mappings has been developed predominantly over the real scalar field, due to the advantages provided by the order structure of $\R$ and the validity of the McShane--Whitney extension theorem, whereas operator space theory has traditionally been primarily formulated over the complex scalars. 
We shall consider both the real and complex settings, which required us to extend to the real case some results about absolute matrix convexity that are usually stated and proved only for complex scalars (see Section \ref{sec-real-absolute-matrix-convexity}). 

In the classical setting, one of the strengths of the Lipschitz-free space is that there are several equivalent approaches to constructing it:
(a) Godefroy and Kalton \cite{Godefroy-Kalton} defined it as the closed linear span of the evaluation functionals, which relies on the Dixmier--Ng theorem;
(b) The approach using de Leeuw’s map realizes the space of Lipschitz functions vanishing at a point as a weak$^*$-closed subspace of an $\ell_\infty$-space, which then yields a representation of the Lipschitz-free space as a quotient of an $\ell_1$-space;
(c) The original Arens--Eells construction \cite{Arens-Eells}, also used in \cite{Weaver-Lipschitz-algebras}, realizes it in terms of molecules, providing a clear and powerful connection to the optimal transportation problem.

In the present setting, we have been able to develop noncommutative counterparts of all these three constructions, providing our Lipschitz-free \emph{operator} spaces with the same rich range of different but equivalent descriptions that has been so fruitful in the classical setting:
(a) In Section \ref{sec-Dixmier-Ng-Kaijser} we identify the Lipschitz-free space as the closure of a linear span of evaluation functionals (Theorem \ref{thm-predual-version-Dixmier-Ng-Kaijser}), using an operator space version of Kaijser's version of the Dixmier--Ng theorem we previously developed in  \cite{chavez2026linearizing};
(b) In Section \ref{sec-de-Leeuw}, we develop an operator space version of de Leeuw's map, which identifies our space of Lipschitz functions as a weak$^*$ closed subspace of an $\ell_\infty$ sum of finite-dimensional matrix algebras. As a consequence, our Lipschitz-free space is described as a quotient of an $\ell_1$-sum of of finite-dimensional Schatten-1 spaces (Proposition \ref{prop-Lip-free-de-Leeuw});
(c) In Section \ref{sec-Arens-Eells}, we describe the Lipschitz free space in terms of molecules as in the Arens--Eells approach. In analogy with the classical case, there is one result describing the operator space structure on the subspace spanned by the molecules (Proposition \ref{prop-Arens-Eells-incomplete-without-matrix} and Corollary \ref{cor-Arens-Eells-incomplete}) and another one for the completion (Theorem \ref{thm-Arens-Eeels}). Just as the classical description has a form reminiscent of the projective tensor product of Banach spaces, in our setting the expressions are related to the projective tensor product of operator spaces.

Finally, in Section \ref{sec-linearization-property} we show that our construction satisfies the desired linearization property (Theorem \ref{thm-linearization}).

Naturally, the present work only lays down the very first step towards a theory of Lispchitz-free spaces in the noncommutative setting.
In a companion paper \cite{CD-Dimant-futuro} we will continue this line of research, concentrating on the case where the distinguished matrix set is the matrix unit ball, and in particular developing a version of the Lipschitz-lifting property of Godefroy and Kalton.

\section{Notation and preliminaries}\label{sec:prelim}

The open unit ball of a Banach space $X$ will be denoted by $B_X$,
and the closed one by $\overline{B}_X$.

\subsection{Operator spaces}
We assume familiarity with the basic theory of operator spaces. The books \cite{Pisier-Operator-Space-Theory} and \cite{Effros-Ruan-book} are excellent references on this topic. While most of the literature deals with complex operator spaces, we will also be considering real ones.
A number of basic properties from the complex case have been known to also hold in the real one for a long time, see \cite{Ruan-real-OS,Ruan-complexifications,Sharma}.
See also \cite{blecher2023real,blecher2024real,blecher2024m,blecher2025realoperatorsystems,blecher2025realnoncommutativeconvexityi,blecher2025realnoncommutativeconvexityii,blecher2025realdecomposablemapsoperator} for a wave of recent interest in the subject.

We denote by $\K$ the scalar field, either $\R$ or $\C$.
A \emph{(concrete) $\K$-operator space} is a closed subspace $V$ of  $\mathcal{B}(H)$, where $H$ is a $\K$-Hilbert space and $\mathcal{B}(H)$ denotes the set of bounded $\K$-linear operators on $H$. Using the identifications $M_n(\mathcal{B}(H)) =\mathcal{B}(H \oplus \cdots \oplus H)$ yields norms on the spaces $M_n(V)$  of $n\times n$ matrices with entries in $V$.
There is also an abstract characterization of operator spaces in terms of properties of these norms, namely Ruan's axioms.

For a function $f : V \to W$ and $n\in\N$, its \emph{$n$-th amplification} is the function $f_n: M_n(V)\to M_n(W)$ given by applying $f$ entry-wise, that is, $f_n\big((x_{jk})_{jk}\big)=(f(x_{jk}))_{jk}$ for all $(x_{jk})_{jk}\in M_n(V)$.
A linear map $T : V \to W$  is    \emph{completely bounded} if 
$\|T\|_{\cb}=\sup_{n\in\N} \|T_n\|<\infty$.
We denote by $\CB(V,W)$ the space of completely bounded mappings from $V$ to $W$, which is an operator space with the structure given by the identification $M_n(\CB(V,W))=\CB(V,M_n(W))$.

The fact that each bounded linear mapping from $V$ to $\K$ is completely bounded provides the dual space $V'$ with the operator space structure given by the identification $M_n(V')=\CB(V,M_n)$. The norm of a matrix in $M_n(V')$ can also be computed through the \textit{matrix pairing} $\langle\langle \cdot, \cdot\rangle\rangle: M_n(V)\times M_n(V')\to M_{n^2}$,  where $\langle\langle (x_{jk}),(\varphi_{rs})\rangle\rangle =(\varphi_{rs}(x_{jk}))$. Precisely,
$$
\|(\varphi_{rs})\|_{M_n(V')}=\sup\big\{\|\langle\langle (x_{jk}),(\varphi_{rs})\rangle\rangle\|_{M_{n^2}}: \|(x_{jk})\|_{M_n(V)}\le 1\big\}.
$$

A mapping $T : V \to W$  is  a  \emph{complete metric surjection} if it is onto and the associated map from $V/\ker(T)$ to $W$ is a completely isometric isomorphism.
Analogously to the Banach case, this is equivalent to having $T_n(B_{M_n(V)}) = B_{M_n(W)}$ for each $n\in\N$.

Given a family $\{V_\gamma\}_{\gamma \in \Gamma}$ of operator spaces, its $\ell_\infty$ direct sum 
 which we will denote by $\ell_\infty( \{V_\gamma \}_{ \gamma \in \Gamma} )$ is the subspace  $V$ of bounded families in   $\bigoplus_{\gamma \in \Gamma} V_\gamma$ with the operator space structure given as follows: 
for $u = (u_\gamma)_{\gamma \in \Gamma} \in  M_n(V)$,
$$
\|u\|_{\ell_\infty( \{V_\gamma \}_{ \gamma \in \Gamma} ) } = \sup \big\{ \|u_\gamma\|_{M_n(V_\gamma)} \;:\; \gamma \in \Gamma \}.
$$

The $\ell_1$-sum of $\{V_\gamma\}_{\gamma \in \Gamma}$ will be similarly denoted by $\ell_1(\{V_\gamma\}_{\gamma \in  \Gamma})$,
and its operator space structure is characterized by the following universal property: for any operator space $W$ and any linear map $T : \ell_1(\{V_\gamma\}_{\gamma \in  \Gamma}) \to W$, $\n{T}_{\cb} \le 1$ if and only if for all $\gamma\in \Gamma$ we have $\n{TJ_\gamma}_{\cb} \le 1$, where $J_\gamma : V_\gamma \to \ell_1(\{V_\zeta\}_{\zeta \in  \Gamma})$ is the canonical injection (that is, the map sending $v \in V_\gamma$ to the vector having $v$ in the $\gamma$-th position and 0 everywhere else).

\subsection{Complexifications}

If $V$ is a real vector space, its \emph{complexification} is the space
\[
V_c = V \oplus i V = \{ x+iy \;:\; x, y \in V \}
\]
with its natural complex linear structure
\begin{align*}
(x_1+iy_1) + (x_2+iy_2) &= (x_1+x_2) + i(y_1+y_2), \qquad x_1,x_2,y_1,y_2 \in V \\
(a + i b) (x+iy) &= (a x- b y) + i (b x + a y), \qquad a,b\in \R, x,y \in V
\end{align*}
and the conjugation map $\conj : V_c \to V_c$
\[
\conj(x+iy) = x - iy, \qquad x,y \in V.
\]
We will denote by $\kappa_V : V \to V_c$ (and almost always just $\kappa$ for simplicity) the canonical inclusion $\kappa_V(x) = x + i0$.

It will sometimes be convenient to observe that the complexification $V_c$ can also be identified with the real tensor product
$V \otimes_{\R} \C$.
Of course, being a real tensor product of real vector spaces, the space $V \otimes_{\R} \C$ is in principle just a real vector space. However, it can be made into a complex vector space by defining
\[
\alpha( v \otimes \beta ) = v \otimes (\alpha\beta) \qquad v\in V, \alpha,\beta \in \C.
\]

When $V$ is additionally a real Banach space, a complex Banach space norm $\n{\cdot}_{V_c}$ on $V_c$ is called a \emph{reasonable complex extension} of the original norm on $V$ if it satisfies for all $x, y \in V$,
\[
\n{x+i0}_{V_c}  = \n{x}, \qquad \n{x+iy}_{V_c} = \n{x-iy}_{V_c}.
\]
There exist a variety of such reasonable complex extensions; we emphasize the \emph{Bochnak complexification} given by
\[
\n{v}_{B} = \inf\bigg\{ \sum_{j=1}^n |\lambda_j| \n{v_j} \; : \; v = \sum_{j=1}^n \lambda_j v_j, \lambda_j \in \C, v_j \in V \bigg\}, \qquad v \in V_c
\]
and the \emph{Taylor complexification} given by
\[
\n{x + iy}_T = \sup_{\theta\in[0,2\pi]} \big\{ \n{(\cos \theta) x + (\sin \theta)y }_V \big\}, \qquad x,y \in V.
\]
Note that under the identification $V_c = V \otimes_\R \C$ mentioned above, the Bochnak norm (resp. the Taylor norm) corresponds to the projective norm $\pi$ (resp. the injective norm $\varepsilon$) on $V \otimes_\R \C$. Once again, in principle the projective and injective tensor norms on $V \otimes_\R \C$ are just real norms, but it is easy to see that in this case they are in fact complex norms with respect to the standard complex structure on $V \otimes_\R \C$.

Yet another useful description of the complexification $V_c$ comes from generalizing the identification of $\C$ with the space of matrices
\[
\left\{ \begin{bmatrix}
a &-b\\
b &a\\
\end{bmatrix} \;:\; a,b\in\R \right\} \subset M_2(\R).
\]
Indeed, $V_c$ can be identified with the space
\[
\left\{ \begin{bmatrix}
x &-y\\
y &x\\
\end{bmatrix} \;:\; x,y\in V \right\} \subset M_2(V)
\]
since whenever $a,b\in\R$, $x, y \in V$ we have
\begin{equation}\label{eqn-eta-and-products}
    \begin{bmatrix}
a &-b\\
b &a\\
\end{bmatrix} \begin{bmatrix}
x &-y\\
y &x\\
\end{bmatrix} = \begin{bmatrix}
ax-by &-(bx+ay)\\
bx+ay &ax-by\\
\end{bmatrix} = \begin{bmatrix}
x &-y\\
y &x\\
\end{bmatrix}\begin{bmatrix}
a &-b\\
b &a\\
\end{bmatrix} .
\end{equation}
Thus, if we denote by $\eta_V : V_c \to M_2(V)$ the map
\[
\eta_V(x + iy) = \begin{bmatrix}
x &-y\\
y &x\\
\end{bmatrix}, \qquad x,y \in V,
\]
we have observed in \eqref{eqn-eta-and-products} that we always have, for $z_1, z_2 \in \C = \R_c$ and any $v \in V_c$,
\[
\eta_V(z_1 v z_2) = \eta_\R(z_1) \eta_V(v) \eta_\R(z_2).
\]
It immediately follows that for any $z_1 \in M_{n,m}(\C)$, $v \in M_{m,k}(V_c)$ and $z_2 \in M_{k,\ell}(\C)$ we have
\begin{equation}\label{eqn-eta-and-matrix-products}
(\eta_V)_{n,\ell}(z_1 \cdot v \cdot z_2) = (\eta_\R)_{n,m}(z_1)  \cdot (\eta_V)_{m,k}(v) \cdot (\eta_\R)_{k,\ell}(z_2).    
\end{equation}

Following \cite{Ruan-complexifications}, a complexification $V_c$ of a real operator space $V$ is \emph{reasonable} if  $\|x+i0\|_n=\|x\|_n$ and $\|x+iy\|_n=\|x-iy\|_n$ for all $x,y\in M_n(V)$ and all $n\in \N$. It is proved in \cite[Thm. 3.1]{Ruan-complexifications} that there is only one reasonable complexification of a real operator space, which is given by the identification $\eta_V$:

\begin{theorem}\label{thm-embedding-complexification-into-M2}
If $V$ is a real operator space, then the map $\eta_V: V_c \to M_2(V)$ 
is a real complete isometry.
\end{theorem}

For a real vector space $V$, we define the real and imaginary part mappings $\Re, \Im : V_c \to V$ by
\[
\Re(x+iy) = x, \quad \Im(x+iy) = y, \qquad x, y \in V.
\]
The $n$-amplifications of these maps, which we denote without the subindex $n$ for simplicity, behave well with real scalars.

\begin{corollary}\label{cor-real-part-is-complete-contraction}
If $V$ is a real operator space, then $\Re, \Im : V_c \to V$
are completely contractive $\R$-linear maps.
\end{corollary}

\begin{proof}
   For the real part, notice that for any $v \in V_c$ we have $\Re(v) = \alpha^* \cdot \eta_V(v) \cdot \alpha$ where $\alpha = \begin{bmatrix}
1\\
0\\
\end{bmatrix}$, showing that $\Re$ is a real linear complete isometry followed by a real linear complete contraction. 
A similar argument works for $\Im(v) = \Re(-i v)$.
\end{proof}

\section{Absolute matrix convexity in the real setting}\label{sec-real-absolute-matrix-convexity}

Absolute matrix convexity in the complex setting was developed in \cite{Effros-Webster}. We will also need the real version, which will be developed in this section.

A \emph{matrix set $\mathbf{A} = (A_n)_n$ over a set $V$} is a sequence of subsets $A_n \subseteq M_n(V)$ for each $n\in\N$.
We sometimes refer to $A_n$ as the $n$-th level of $\mathbf{A}$.
For a map $T : V \to W$, the expression $T(\mathbf{A})$ denotes the matrix set over $W$ given by $\big(T_n(A_n)\big)_n$  where $T_n$ is the $n$-th amplification of $T$.
Given two matrix sets $\mathbf{A} = (A_n)_n$ and $\mathbf{C} = (C_n)_n$ over the same set $V$, we write $\mathbf{A} \subseteq \mathbf{C}$ when $A_n \subseteq C_n$ for all $n\in\N$.
When $V$ is an operator space, we say that  $\mathbf{A}$ is open (resp. closed) whenever each $A_n$ is open (resp. closed). By the closure of a matrix set $\mathbf{A} = (A_n)_n$ we mean taking the closure on each level: $\overline{\mathbf{A}} = ( \overline{A_n})_n$.
We say that $\mathbf{A}$ is \emph{bounded} if there exists a constant $C$ such that for all $n\in \N$ and all $x \in A_n$ we have $\n{x}_{M_n(V)} \le C$.

Sometimes we will consider a finite sequence $\mathbf{A} = (A_n)_{n=1}^N$ with $A_n \subseteq M_n(V)$ as a matrix set, by implicitly defining $A_n = \emptyset$ whenever $n>N$.
In particular, a subset $A_1$ of $V$ will be considered as a matrix set over $V$
 by taking $A_n = \emptyset$ for $n>1$.

We say that a matrix set $\mathbf{A} = (A_n)$ over a $\K$-vector space $V$ is \emph{$\K$-absolutely matrix convex} if:
\begin{enumerate}[(i)]
	\item For all $x \in A_n$ and $y \in A_m$, $x \oplus y \in A_{m+n}$.
	\item For all $x \in A_n$, $\alpha \in M_{m,n}(\K)$, and $\beta \in M_{n,m}(\K)$ with $\n{\alpha}, \n{\beta} \le 1$,
	$\alpha x \beta \in A_m$.
\end{enumerate}

The main example of an absolutely matrix convex set of matrices over a $\K$-operator space $V$ comes from the operator space structure: it follows immediately from Ruan's axioms that
both the \emph{open matrix unit ball of $V$} given by
$\mathbf{B}_V = \big( B_{M_n(V)} \big)_n$
and the 
\emph{closed matrix unit ball of $V$} given by
$\overline{\mathbf{B}}_V = \big( \overline{B}_{M_n(V)} \big)_n$
are $\K$-absolutely matrix convex matrix sets over $V$.

To simplify some statements, whenever $\mathbf{A} = (A_n)_n$ is a matrix set we will write ``$x \in \mathbf{A}$'' as a shorthand for ``there exists $n\in\N$ such that $x\in A_n$''.

Clearly, the intersection of a family of $\K$-absolutely matrix convex matrix sets over $V$ (understood as intersecting on each matrix level) is again absolutely matrix convex. Therefore we can define the \emph{$\K$-absolutely matrix convex hull} of a given matrix set $\mathbf{A}$ over $V$, denoted by $\amconv_\K(\mathbf{A})$,
as the smallest $\K$-absolutely matrix convex matrix set over $V$ that contains $\mathbf{A}$.
In order to better understand the next result, recall that if $A$ is a subset of a $\K$-vector space $V$, then the absolutely convex hull of $A$ can be described as
\[
\aconv(A)=\bigg\{ \sum_{j=1}^m \lambda_j x_ j \; : \; x_j \in A, \lambda_j \in \K,  \sum_{j=1}^m |\lambda_j| \le 1 \bigg\}.
\]
Furthermore, the Cauchy--Schwarz inequality implies that an equivalent description is
\begin{equation}\label{eqn-abs-convex-hull}
    \aconv(A)=\bigg\{ \sum_{j=1}^m \alpha_jx_ j\beta_j \; : \; x_j \in A, \alpha_j,\beta_j \in \K,  \sum_{j=1}^m |\alpha_j|^2 \le 1, \sum_{j=1}^m |\beta_j|^2 \le 1 \bigg\}.
\end{equation}
The following nice characterization of the absolutely matrix convex hull of a matrix set, clearly analogous to the above expression, is due to B.E. Johnson \cite[Lemma 3.2]{Effros-Webster} (the slight restatement below is from \cite[Lemma 2.1]{CD-Oikhberg}).
We point out that while \cite{Effros-Webster} assumes complex scalars, the same proof works in the real case.

\begin{lemma}\label{lemma-characterization-abs-mat-convex-hull}
	Given a matrix set $\mathbf{A}=(A_n)_n$ over $V$, the $n$-th level of $\amconv_\K(\mathbf{A})$ is given by
\[
\bigg\{ \sum_{j=1}^m \alpha_j x_j \beta_j \; : \; x_j \in A_{k_j}, \alpha_j \in M_{n,k_j}(\K),\beta_j \in M_{k_j,n}(\K) \text{ such that } \bigg\| \sum_{j=1}^m \alpha_j\alpha_j^*  \bigg\| \le 1,  \bigg\|\sum_{j=1}^m \beta_j^*\beta_j \bigg\| \le 1  \bigg\}.
\]
\end{lemma}

We remark that the sum $\sum_{j=1}^m \alpha_jx_ j\beta_j$ appearing in  \eqref{eqn-abs-convex-hull} can be rewritten 
as a matrix product of the form $\alpha \cdot D \cdot \beta$: set
\begin{align*}
\alpha &:= \begin{bmatrix}
\alpha_1 & \alpha_2 &\cdots &\alpha_m\\
\end{bmatrix} \in M_{n,k}(\K), \quad
\beta = \begin{bmatrix}
\beta_1 \\ \beta_2 \\ \vdots \\ \beta_m\\
\end{bmatrix} \in M_{k,n}(\K), \\ 
D &:= \diag(x_1, x_2, \dotsc, x_m) \in M_k(V),
\end{align*}
where $k=k_1+k_2+\cdots+k_m$.
Additionally, the conditions $\bigg\| \sum_{j=1}^m \alpha_j\alpha_j^*  \bigg\| \le 1,  \bigg\|\sum_{j=1}^m \beta_j^*\beta_j \bigg\| \le 1 $ are equivalent to $\n{\alpha} \le 1$, $\n{\beta} \le 1$.
Furthermore, any representation as a matrix product $\alpha \cdot D \cdot \beta$ of this form corresponds to a sum $\sum_{j=1}^m \alpha_j x_j \beta_j$.
Therefore, Lemma \ref{lemma-characterization-abs-mat-convex-hull}, admits the following equivalent statement.

\begin{corollary}\label{cor-amconvex hull of increasing matrix form}
Given a matrix set $\mathbf{A}=(A_n)_n$ over $V$, the $n$-th level of $\amconv_\K(\mathbf{A})$ is given by the products of the form $\alpha \cdot D \cdot \beta$
where $\alpha \in M_{n,k}(\K)$ and $\beta \in M_{k,n}(\K)$ satisfy $\n{\alpha} \le 1$, $\n{\beta} \le  1$, and $D \in M_k(V)$ is a block diagonal matrix with diagonal entries $x_j \in A_{k_j}$ for $1 \le j \le m$, and $k=k_1+k_2+\cdots+k_m$.
\end{corollary}

\begin{definition}
A matrix set $\mathbf{A}$ over a complexification of a real space is said to be \emph{symmetric} if $\conj(\mathbf{A}) = \mathbf{A}$, where $\conj$ is the conjugation. 
\end{definition}

\begin{lemma}\label{lemma-hull-of-symmetric}
If a matrix set $\mathbf{A}$ over a complexified operator space is symmetric, then $\amconv_\C(\mathbf{A})$ is also symmetric.
\end{lemma}

\begin{proof}
Let $x \in \amconv_\C(\mathbf{A})$. By Lemma \ref{lemma-characterization-abs-mat-convex-hull}, there exist  $x_j \in A_{k_j}$, $\alpha_j \in M_{n,k_j}(\C)$, $\beta_j \in M_{k_j,n}(\C)$ such that $\big\| \sum_{j=1}^m \alpha_j\alpha_j^*  \big\| \le 1$, $\big\|\sum_{j=1}^m \beta_j^*\beta_j \big\| \le 1 $ and $x = \sum_{j=1}^m \alpha_j x_j \beta_j$.
Conjugating we get $\conj(x) = \sum_{j=1}^m \conj(\alpha_j)\, \conj(x_j) \conj(\beta_j)$, from where it is clear that $\conj(x) \in \amconv_\C(\mathbf{A})$ since $\conj(x_j) \in \mathbf{A}$.
\end{proof}

We now state a matrix analogue of the fact that taking the convex hull of a set cannot increase the supremum of the norm of the vectors in the set.
For a matrix set $\mathbf{A}=(A_n)_n$ over an operator space $V$, we will denote
$\n{\mathbf{A}} = \sup\{ \n{x}_{M_n(V)} \;:\; n\in\N, x \in A_n \}$.
It was proved in \cite[Lemma 2.2]{chavez2026linearizing} for complex scalars, but the same proof works in the real case.

\begin{lemma}\label{lemma-norm-of-abs-mat-conv-hull}
For a matrix set $\mathbf{A}=(A_n)_n$ over a $\K$-operator space $V$, $\n{\mathbf{A}} = \n{\amconv_\K(\mathbf{A})}$.
\end{lemma}

The following two results relate the real and complex absolutely matrix convex hull procedures in a real operator space and its complexification.
The first one basically says that if we start with real objects, then the real part of a complex absolutely matrix convex combination is in fact a real absolutely matrix convex combination.

\begin{lemma}\label{lemma-real-part-of-complex-convex-hull}
 Let $V$ be a real vector space.
 For any matrix set $\mathbf{A}$ over $V$ we have 
 \[
 \Re\big(  \amconv_\C(\kappa(\mathbf{A}))  \big) \ \subseteq \amconv_\R(\mathbf{A}).
 \]
\end{lemma}

\begin{proof}
Consider an element $x \in M_n(V_c)$ belonging to the $n$-th level of $\amconv_\C(\kappa(\mathbf{A}))$.
By Lemma \ref{lemma-characterization-abs-mat-convex-hull}, there exists a representation  $x = \sum_{j=1}^m \alpha_j  x_j  \beta_j$ where  $x_j \in A_{k_j}$ (note we are, in the usual abuse of notation, identifying $x_j \in M_{k_j}(V)$ with $\kappa_{k_j}(x_j) \in M_{k_j}(V_c)$), $\alpha_j \in M_{n,k_j}(\C)$, $\beta_j \in M_{k_j,n}(\C)$ satisfy $\bigg\| \sum_{j=1}^m \alpha_j\alpha_j^*  \bigg\| \le 1$,  $\bigg\|\sum_{j=1}^m \beta_j^*\beta_j \bigg\| \le 1$.
It now follows that
\[
x = \sum_{j=1}^m \big( \Re(\alpha_j) x_j \Re(\beta_j)- \Im(\alpha_j)x_j \Im(\beta_j)\big) + i \sum_{j=1}^m \big( \Im(\alpha_j) x_j \Re(\beta_j)+ \Re(\alpha_j)x_j \Im(\beta_j)\big).
\]
We then conclude
\begin{equation}\label{eqn-relation-convex-hulls-real-part}
\Re(x) = \sum_{j=1}^m \big( \Re(\alpha_j) x_j \Re(\beta_j)- \Im(\alpha_j)x_j \Im(\beta_j)\big)
\end{equation}
Note that
\begin{align*}
   \sum_{j=1}^m \alpha_j\alpha_j^* &= \sum_{j=1}^m \big( \Re(\alpha_j)+i\Im(\alpha_j) \big)\big( \Re(\alpha_j)^*-i\Im(\alpha_j)^* \big) \\
&= \sum_{j=1}^m \big( \Re(\alpha_j)\Re(\alpha_j)^* + \Im(\alpha_j)\Im(\alpha_j)^* \big) + i \sum_{j=1}^m \big( -\Re(\alpha_j)\Im(\alpha_j)^* + \Im(\alpha_j)\Re(\alpha_j)^* \big),
\end{align*}
so that by Corollary \ref{cor-real-part-is-complete-contraction} we get
\begin{equation}\label{eqn-relation-convex-hulls-alpha}
    \n{ \sum_{j=1}^m \big( \Re(\alpha_j)\Re(\alpha_j)^* + \Im(\alpha_j)\Im(\alpha_j)^* \big) } = \n{\Re\bigg(\sum_{j=1}^m \alpha_j\alpha_j^*\bigg)} \le \n{ \sum_{j=1}^m \alpha_j\alpha_j^* } \le 1.
\end{equation}
Analogously,
\begin{equation}\label{eqn-relation-convex-hulls-beta}
    \n{ \sum_{j=1}^m \big( \Re(\beta_j)^*\Re(\beta_j) + \Im(\beta_j)^*\Im(\beta_j) \big) } = \n{\Re\bigg( \sum_{j=1}^m \beta_j^*\beta_j \bigg)} \le \n{ \sum_{j=1}^m \beta_j^*\beta_j } \le 1
\end{equation}

It now follows from
\eqref{eqn-relation-convex-hulls-real-part}
\eqref{eqn-relation-convex-hulls-alpha}
and
\eqref{eqn-relation-convex-hulls-beta}, together with 
Lemma \ref{lemma-characterization-abs-mat-convex-hull}, that $\Re(x) \in \amconv_\R( \mathbf{A} )$.    
\end{proof}

The following proposition says that while the absolutely matrix convex hull in the complexification will be larger, it will not pick up any more ``real" parts.

\begin{proposition}\label{prop-complex-hull-of-real-matrix-set}
Let $V$ be a real vector space. For any matrix set $\mathbf{A}$ over $V$ we have
\[
\kappa\big(\amconv_\R(\mathbf{A})\big) = \big\{ x \in \amconv_\C(\kappa(\mathbf{A})) \; : \; \Im(x)=0 \big\}.
 \]
If $V$ is an operator space, we furthermore have
\[
\kappa\big(\overline{\amconv_\R(\mathbf{A})}\big) = \big\{ x \in \overline{ \amconv_\C(\kappa(\mathbf{A}))} \; : \; \Im(x)=0 \big\}.
 \]
\end{proposition}

\begin{proof}
It follows immediately from Lemma \ref{lemma-characterization-abs-mat-convex-hull} that $\kappa(\amconv_\R(\mathbf{A})) \subseteq \amconv_\C(\kappa(\mathbf{A}))$, which gives the inclusion $\subseteq$.

For the containment $\supseteq$, consider an element $x$ of the $n$-th level of $\amconv_\C(\kappa(\mathbf{A}))$ with $\Im(x)=0$.
By Lemma \ref{lemma-real-part-of-complex-convex-hull} we have $\Re(x) \in \amconv_\R(\mathbf{A})$, and hence $x \in \kappa\big( \amconv_\R\mathbf{A} \big)$.

The statement concerning closures follows first in one direction from the continuity of $\kappa$. To obtain the reverse inclusion, let $x_0$ belong to the $n$-th level of  $\overline{\amconv_\C(\kappa(\mathbf{A}))}$ with $\Im(x_0)=0$.
Then, there exists a sequence $(x_j)_{j=1}^\infty$ in the $n$-th level of $\amconv_\C(\kappa(\mathbf{A}))$ converging to $x_0$. Since taking real parts is a complete contraction by Corollary \ref{cor-real-part-is-complete-contraction}, $\Re(x_j)$ converges to $\Re(x_0)$.
Moreover, $\kappa(\mathbf{A})$ is symmetric, and hence so is $\amconv_\C(\kappa(\mathbf{A}))$ by Lemma \ref{lemma-hull-of-symmetric}. Consequently, $\kappa_n(\Re(x_j)) = \frac{1}{2}( x_j + \conj(x_j) ) \in \amconv_\C(\kappa(\mathbf{A}))$, and its imaginary part is zero.
By the first part,  $\Re(x_j) \in \amconv_\R(\mathbf{A})$ for every $j$ and so, $\Re(x_0) \in \overline{\amconv_\R(\mathbf{A})}$, which then satisfies $x_0=\kappa_n(\Re(x_0)) \in \kappa\big(\overline{\amconv_\R(\mathbf{A})}\big)$, finishing the proof.
\end{proof}

We now state two versions of the Hahn--Banach separation theorem for absolute matrix convexity.
In the complex case both are known: Theorem \ref{thm-Hahn-Banach-matrix-convexity} is stated in \cite[Thm. 2.3]{CD-Oikhberg} and is contained in the proof of \cite[Prop. 4.1]{Effros-Webster}, while Theorem \ref{thm-Hahn-Banach-operator-TVS} is \cite[Lemma 3.1.8]{ChaDiGa-tensor-norms}.
We remark that the same proofs work in the real case.
Alternatively,
one could deduce the real versions from the complex ones by passing to the complexification, but this produces maps with values in $M_n(\C) \equiv M_{2n}(\R)$, thus doubling the size of the matrices.

\begin{theorem}\label{thm-Hahn-Banach-matrix-convexity}
	Let $\mathbf{A} = (A_n)_n$ be a closed absolutely matrix convex  matrix set over a $\K$-operator space $V$ and let $x_0 \in M_n(V) \setminus A_n$ for some $n\in\N$.
	Then there exist $k\in\N$ and $\varphi \in \CB(V,M_k(\K))$  such that for all $m\in\N$ and all $x \in A_m$,
	$$
	\n{ \varphi_m(x) }_{M_{mk}(\K)} \le 1 \quad \text{but} \quad \n{ \varphi_n(x_0) }_{M_{nk}(\K)} > 1.
	$$
\end{theorem}

\begin{theorem}\label{thm-Hahn-Banach-operator-TVS}
Let $V$ be a $\K$-vector space, and $\mathbf{A}=(A_n)_{n=1}^\infty$ an absolutely matrix convex matrix set over $V$ such that $A_1$ is absorbent. Suppose that there exist $u \in M_n(V)$ and $r>0$ such that $(u +rA_n) \cap A_n = \varnothing$. Then there exist $k\in\N$ and a linear map $\varphi : V \to M_k(\K)$ such that $\n{\varphi_n(u)}_{M_{nk}(\K)} >1$, and $\n{\varphi_m(v)}_{M_{mk}(\K)} \le 1$ for every $m\in \N$ and every $v \in A_m$.     
\end{theorem}

For a real Banach space $V$, we denote by $V_{c,B}$ its Bochnak complexification, given by, for $z \in V_c$,
\[
\n{z}_{B} = \inf\bigg\{ \sum_{j=1}^n |\lambda_j| \n{x_j} \; : \; z = \sum_{j=1}^n \lambda_j x_j, \lambda_j \in \C, x_j \in V \bigg\}.
\]

The following proposition is simply a restatement of the definition above.

\begin{proposition}\label{prop-ball-Bochnak}
For a real Banach space $V$ we have
$B_{V_{c,B}} = \aconv_\C(\kappa(B_V))$.
\end{proposition}

We now prove an analogue of the previous result, in the operator space setting.
It interprets the canonical complexification of a real operator space in terms of a convex hull procedure.

\begin{theorem}\label{thm-complexification-as-hull}
Let $V$ be a real operator space and $V_c$ its canonical reasonable complexification. Then
\[
\mathbf{B}_{V_c} = \amconv_\C (\kappa( \mathbf{B}_V)).
\]
\end{theorem}

\begin{proof}
It is clear that $\kappa( \mathbf{B}_V) \subseteq \mathbf{B}_{V_c}$, and therefore  $\amconv_\C(\kappa( \mathbf{B}_V)) \subseteq \mathbf{B}_{V_c}$.
Consider the complex matrices $\alpha = \begin{bmatrix}
1/\sqrt{2} & i/\sqrt{2}\\
\end{bmatrix}$, $\beta = \begin{bmatrix}
1/\sqrt{2} \\ -i/\sqrt{2}
\end{bmatrix}$, observe that $\n{\alpha}, \n{\beta} = 1$.
Note that for any $x,y \in V$,
\[
\alpha \begin{bmatrix}
x &-y\\
y &x\\
\end{bmatrix} \beta = x + iy.
\]
Therefore, recalling that $\eta_V : V_c \to M_2(V)$ given by
\[
\eta_V(x + iy) =  \begin{bmatrix}
x &-y\\
y &x\\
\end{bmatrix},
\]
is a completely isometric embedding (Theorem \ref{thm-embedding-complexification-into-M2})
we see that $B_{V_c} \subseteq \amconv_\C (\kappa_2( B_{M_2(V)} )) \subseteq \amconv_\C (\kappa( \mathbf{B}_{V} ))$.
Similarly, for a matrix $ x \in M_n(V_c)$, we have
$x = (\alpha \otimes Id) (\eta_V)_n(x) (\beta \otimes Id)$. Since $\n{\alpha \otimes Id} = \n{\alpha} \n{Id} = 1$ and $\n{\beta \otimes Id} = \n{\beta} \n{Id} = 1$, we then conclude that $B_{M_n(V_c)} \subseteq \amconv_\C (\kappa_{2n}( B_{M_{2n}(V)} )) \subseteq \amconv_\C (\kappa( \mathbf{B}_{V} ))$, which implies $\mathbf{B}_{V_c} \subseteq \amconv_\C \kappa( \mathbf{B}_V)$.
\end{proof}

Recall that an operator space $V$ is said to be \emph{maximal} if for any operator space $W$ and any linear map $T : V \to W$ we have $\n{T : V \to W}_{\cb} = \n{T : V \to W}$.
Given a Banach space $V$, there is a unique maximal operator space structure on $V$, denoted $\MAX(V)$, which is compatible with the norm on $V$.
The following characterization of the matrix unit ball of the maximal operator space structure associated to a Banach space is well known to specialists. The complex version can be found as a particular case of \cite[Lemma 8.10]{ChaDiGa-Revisiting}. Recall that $\amconv_\K\big( B_V \big)$ means the $\K$-absolutely matrix convex hull of the matrix set $\mathbf{A}=(A_n)_{n=1}^\infty$ with $A_1=B_V$ and $A_n=\emptyset$ for all $n>1$. 

\begin{lemma}\label{lemma-matrix-ball-of-maximal}
    Let $V$ be a $\K$-Banach space. Then $\overline{\mathbf{B}}_{\MAX(V)} = \overline{\amconv_\K\big( B_V \big)}$.
\end{lemma}

\begin{proof}
Since $\overline{\mathbf{B}}_{\MAX(V)} \supseteq \overline{B}_V$, we immediately conclude the containment $\supseteq$.
If the containment were strict, Theorem \ref{thm-Hahn-Banach-matrix-convexity} yields a contradiction.
\end{proof}

The following is a version of Proposition \ref{prop-ball-Bochnak} for maximal operator spaces.

\begin{corollary} \label{cor-MAX-Bochnak}
For a real Banach space $V$,
$\MAX( V_{c,B} ) = \MAX(V)_c$ completely isometrically.   
\end{corollary}

\begin{proof}
Using Proposition \ref{prop-ball-Bochnak} and the fact that $\kappa$ is $\R$-linear,
\begin{align*}
\amconv_\C (B_{V_{c,B}}) &=  \amconv_\C ( \aconv_\C(\kappa(B_V))) 
=  \amconv_\C  (\kappa(B_V)) \\
&=  \amconv_\C (\amconv_\R (\kappa( B_V))) = 
 \amconv_\C  (\kappa(\amconv_\R (B_V))).
\end{align*}
Taking closures and using Lemma \ref{lemma-matrix-ball-of-maximal} we thus get
\[
\overline{\mathbf{B}}_{\MAX( V_{c,B} )} = \overline{ \amconv_\C  (\kappa(\amconv_\R (B_V))) }.
\]
On the other hand,
by Theorem \ref{thm-complexification-as-hull} we know 
\[
\overline{\mathbf{B}}_{\MAX(V)_c} = \overline{ \amconv_\C  (\kappa(\mathbf{B}_{\MAX(V)}))  } =   \overline{ \amconv_\C  (\kappa(\overline{\mathbf{B}}_{\MAX(V)}))  },
\]
and using Lemma \ref{lemma-matrix-ball-of-maximal} again thus yields
\[
\overline{\mathbf{B}}_{\MAX(V)_c} =  \overline{ \amconv_\C  (\kappa( \overline{ \amconv_\R (B_V) } ))  } = \overline{ \amconv_\C  (\kappa( \amconv_\R (B_V) ))  }.
\]
We conclude $\overline{\mathbf{B}}_{\MAX( V_{c,B} )} = \overline{\mathbf{B}}_{\MAX(V)_c}$, finishing the proof.
\end{proof}

More generally, an operator space $V$ is said to be \emph{$N$-maximal} if for any operator space $W$ and any linear map $T : V \to W$ we have $\n{T : V \to W}_{\cb} = \n{T_N : M_N(V) \to M_N(W)}$.
Once again, the following is well known.
A version for complex scalars has appeared in \cite[Lemma 8.10]{ChaDiGa-Revisiting}, and the same argument works in the real case.

\begin{lemma}\label{lemma-matrix-ball-of-N-maximal}
    Let $V$ be a $\K$-operator space. Then $V$ is $N$-maximal if and only if
    \[
    \overline{\mathbf{B}}_{V} = \overline{\amconv_\K\big( B_{M_n(V)} \big)_{n=1}^N}
    \]
\end{lemma}

\begin{lemma}\label{lemma-convex-hull-of-finite-height}
 Suppose that for a $\K$-operator space $V$ we have  $\overline{\mathbf{B}}_V = \overline{\amconv_\K \mathbf{A}}$, where $\mathbf{A}$ is a matrix set over $V$ such that $A_n = \emptyset$ for $n>N$. Then $V$ is $N$-maximal.
\end{lemma}

\begin{proof}
For any operator space $W$ and any continuous linear map $T : V \to W$ we have by the linearity of $T$ and Lemma \ref{lemma-norm-of-abs-mat-conv-hull},
\begin{align*}
    \n{T}_{\cb} &= \n{ T( \overline{\mathbf{B}}_V ) } = \n{ T \left(\overline{\amconv_\K \big( A_n \big)_{n=1}^N }\right) } 
    = \n{ T \left(\amconv_\K \big( A_n \big)_{n=1}^N \right) } 
    = \n{ \amconv_\K T \big(A_n \big)_{n=1}^N } \\
    &= \n{ \amconv_\K\big( T_nA_n \big)_{n=1}^N } \le \n{T_N : M_N(V) \to M_N(W)} \le \n{T}_{\cb},
\end{align*}
showing that $V$ is $N$-maximal.
\end{proof}

\begin{corollary}\label{cor: N maximal complexification}
If a real operator space $V$ is $N$-maximal, then $V_c$ is an $N$-maximal complex operator space.    
\end{corollary}

\begin{proof}
By Theorem \ref{thm-complexification-as-hull} and Lemma \ref{lemma-matrix-ball-of-N-maximal}, together with the fact that $\kappa$ is $\R$-linear and continuous,
\begin{align*}
\overline{\mathbf{B}}_{V_c} &= \overline{\amconv_\C (\kappa(\overline{\mathbf{B}}_V)) }
= \overline{\amconv_\C \left(\kappa \Big( \overline{\amconv_\R \big( B_{M_n(V)} \big)_{n=1}^N}\Big)\right) }
= \overline{\amconv_\C \left(\kappa \Big(\amconv_\R \big( B_{M_n(V)} \big)_{n=1}^N \Big)\right)} \\
&= \overline{\amconv_\C \left( \amconv_\R \Big(\kappa \big( B_{M_n(V)} \big)_{n=1}^N \Big)\right)}
= \overline{\amconv_\C    \big( \kappa_n (B_{M_n(V)}) \big)_{n=1}^N }
\end{align*}
which implies that $V_c$ is $N$-maximal by Lemma \ref{lemma-convex-hull-of-finite-height}.
\end{proof}

\begin{remark}
   There is a large body of results where a real operator space $V$ has a certain property if and only if $V_c$ has the complex version of said property, see e.g. \cite[Sec. 4]{Ruan-complexifications} for examples such as injectivity, nuclearity, exactness, and various approximation properties; see also \cite[Sec. 5]{blecher2023real}, \cite[Sec. 4]{blecher2024real}.
For $N$-maximality the situation is different: the converse of Corollary \ref{cor: N maximal complexification} holds only up to a factor of $2$ in the height, and this factor cannot be removed in general (at least for $N=1$).
This is similar to \cite[Cor. 4.2]{MR4659878}, where a factor of 2 in the size of the amplifications also appears when connecting properties of an object and its complexification.
\end{remark}

\begin{proposition}\label{prop-Vc-N-maximal-implies-V-2N-maximal}
Let $V$ be a real operator space and $N\in\N$. If $V_c$ is $N$-maximal, then $V$ is $2N$-maximal.
\end{proposition}

\begin{proof}
Let $W$ be a real operator space and take $S : V\to W$ any $\R$-linear map. Consider its complexification $S_c : V_c\to W_c$ given by $S_c(x+iy) = Sx+iSy$, which is $\C$-linear.
By \cite[Thm. 2.1]{Ruan-complexifications} we know that $\n{S}_{\cb} = \n{S_c}_{\cb}$ and so  the $N$-maximality of $V_c$ implies, using also that $\eta$ isometrically identifies $M_N(V_c)$ with a subspace of $M_{2N}(V)$,
\[
\n{S}_{\cb} = \n{S_c}_{\cb} = \n{(S_c)_N : M_N(V_c) \to M_N(W_c)} \le \n{S_{2N} : M_{2N}(V) \to M_{2N}(W)} \le \n{S}_{\cb}.
\]
Hence $\n{S}_{\cb} = \n{S_{2N}}$, proving that $V$ is $2N$-maximal.
\end{proof}

The following construction shows that the previous proposition is sharp for $N=1$.
Recall that for a complex Banach space $E$, the matrix norms of $\MAX(E)$ are given, for $u\in M_n(E)$, by (the supremum being over all complex Hilbert spaces $H$)
\[
\n{u}_{M_n(\MAX(E))} = \sup\big\{ \n{T_n(u)}_{M_n(\B(H))} \;:\; T : E\to\B(H) \text{ linear}, \ \n{T}\le 1 \big\}.
\]
Any complex operator space is in particular a real operator space (with the same matrix norms), and so is any real subspace of it.

\begin{proposition}\label{prop-converse-N-maximal-fails}
Let $X$ be a real Banach space, and let $E$ denote $X_c$ endowed with a reasonable complex extension of the norm of $X$ which is different from the Bochnak norm.
Let $V$ be the real operator space obtained by endowing $X$ with the structure inherited from $\MAX(E)$ via $\kappa$, that is,
\[
\n{x}_{M_n(V)} = \n{\kappa_n(x)}_{M_n(\MAX(E))}, \qquad x\in M_n(X).
\]
Then:
\begin{enumerate}[(a)]
\item $V_c = \MAX(E)$ completely isometrically; in particular, $V_c$ is maximal.
\item $V$ is not maximal.
\item $V$ is $2$-maximal.
\end{enumerate}
\end{proposition}

\begin{proof}
$(a)$
By construction, $\MAX(E)$ is a complex operator space structure on the complex vector space $X_c$ whose restriction to $M_n(X)$ (via $\kappa_n$) is the given structure on $M_n(V)$.
By \cite[Thm. 3.1]{Ruan-complexifications} it suffices to check that it is reasonable, i.e., 
 that $\conj : \MAX(E) \to \MAX(E)$ is completely isometric (since $\conj$ is not complex linear, this does not follow immediately from maximality).
Since $\conj$ is its own inverse, it suffices to show that $\conj$ is completely contractive.
By Lemma \ref{lemma-matrix-ball-of-maximal}, we need to show
\[
\conj\Big(\overline{\amconv_\C\big( B_E \big)}\Big) \subseteq \overline{\amconv_\C\big( B_E \big)}.
\]
Since $\conj$ is continuous, it will suffice to show
\[
\conj\Big({\amconv_\C\big( B_E \big)}\Big) \subseteq \overline{\amconv_\C\big( B_E \big)}.
\]
Now, by Lemma \ref{lemma-characterization-abs-mat-convex-hull} an element $x$ of the $n$-th level of $\amconv_\C\big( B_E \big)$ has the form $\sum_{j=1}^m \alpha_j x_j \beta_j$ for some $\alpha_j \in M_{n,1}(\C)$, $\beta_j \in M_{1,n}(\C)$ and $x_j \in B_E$.
Thus $\conj_n(x) = \sum_{j=1}^m \conj(\alpha_j) \conj(x_j) \conj(\beta_j)$, which is clearly in $\amconv_\C\big( B_E \big)$ by Lemma \ref{lemma-characterization-abs-mat-convex-hull} again.

$(b)$
Since the norm of $E$ is a reasonable complex extension of the norm of $X$, the first level of $V$ is $X$ itself.
If $V$ were maximal, then $V=\MAX(X)$ completely isometrically, and Corollary \ref{cor-MAX-Bochnak} would give $V_c = \MAX(X)_c = \MAX(X_{c,B})$.
Looking at the first level and using (a), the norm of $E$ would coincide with the Bochnak norm, contradicting our assumption.

$(c)$
This follows from $(a)$ and Proposition \ref{prop-Vc-N-maximal-implies-V-2N-maximal} with $N=1$.
\end{proof}

For a concrete example, let $X=\ell_1^2$ (real) and let $E$ be its Taylor complexification.
For $z = e_1+ie_2$ we have $\n{z}_T = \sup_\theta \big(|\cos\theta|+|\sin\theta|\big) = \sqrt{2}$.
On the other hand, the norm of $\ell_1^2(\C)$ is a reasonable complex extension of the norm of $X$, and whenever $z=\sum_j\lambda_jx_j$ with $\lambda_j\in\C$, $x_j\in X$, we have $\n{z}_{\ell_1^2(\C)} \le \sum_j|\lambda_j|\n{x_j}$; hence $\n{z}_B\ge\n{z}_{\ell_1^2(\C)} = 2$.
Thus the Taylor and Bochnak norms differ, and Proposition \ref{prop-converse-N-maximal-fails} applies.

\begin{remark}
Proposition \ref{prop-converse-N-maximal-fails} can be summarized as follows: the complexification of a maximal real operator space is maximal and has the Bochnak norm on its first level (Corollary \ref{cor-MAX-Bochnak}); a complex maximal operator space whose first level is any other reasonable complexification is still the complexification of a real operator space, but the latter fails to be maximal.
\end{remark}

It is  well known that the closed unit ball of $\ell_1$ is the closure of the absolutely convex hull of the canonical basis. More generally, the closed unit ball of an $\ell_1$-sum of Banach spaces $\ell_1( \{ X_\gamma\}_{\gamma \in\Gamma} )$ is the closure of the absolutely convex hull of the union of the individual balls $B_{X_\gamma}$ (more precisely, the union of the unit balls of the canonical copies of the $X_\gamma$ inside $\ell_1( \{ X_\gamma\}_{\gamma \in\Gamma} )$).
The following is an operator space version of that fact.

\begin{proposition}\label{prop-matrix-unit-ball-of-ell1-sum}
 Let $\{V_\gamma\}_{\gamma \in\Gamma}$ be a collection of operator spaces. For each $\gamma_0\in \Gamma$, let $J_{\gamma_0} : V_{\gamma_0} \to \ell_1( \{ V_\gamma\}_{\gamma \in \Gamma} )$ be the canonical complete isometry, that is, the one that sends $v \in V_{\gamma_0}$ to the vector with $v$ in the $\gamma_0$-th position and 0 everywhere else.
Then the closed matrix unit ball of $\ell_1( \{ V_\gamma\}_{\gamma \in \Gamma} )$ is the closure of the absolutely matrix convex hull of $\bigcup_{\gamma\in\Gamma} J_\gamma(\mathbf{B}_{V_\gamma})$.
\end{proposition}

\begin{proof}
     The complex version appeared in \cite[Prop. 8.5]{ChaDiGa-Revisiting}. The same proof, based on a Hahn--Banach separation argument, works in the real case thanks to Theorem \ref{thm-Hahn-Banach-matrix-convexity}.
\end{proof}

The next corollary was proved by different methods in \cite[Lemma 5.3]{blecher2023real}.

\begin{corollary}
    For any collection $\{V_\gamma\}_{\gamma \in\Gamma}$ of real operator spaces,
    \[
    \bigg( \ell_1\big(\{V_\gamma\}_{\gamma \in\Gamma}) \bigg)_c = \ell_1\big(\{(V_\gamma)_c\}_{\gamma \in\Gamma} \big) 
    \]
\end{corollary}

\begin{proof}

Using Theorem \ref{thm-complexification-as-hull} and Proposition \ref{prop-matrix-unit-ball-of-ell1-sum}, each of them twice,
\begin{align*}
\overline{\mathbf{B}}_{\big( \ell_1\big(\{V_\gamma\}_{\gamma \in\Gamma}\big) \big)_c} 
&= \overline{ \amconv_\C \left(\kappa( \mathbf{B}_{\ell_1\big(\{V_\gamma\}_{\gamma \in\Gamma}\big)})\right) } = \overline{ \amconv_\C \left(\kappa( \overline{\mathbf{B}}_{\ell_1\big(\{V_\gamma\}_{\gamma \in\Gamma}\big)})\right) } \\
&= \overline{ \amconv_\C \left(\kappa( \overline{ \amconv_\R \bigcup_{\gamma\in\Gamma} J_\gamma(\mathbf{B}_{V_\gamma}) }) \right)} = \overline{ \amconv_\C \left(\kappa(  \bigcup_{\gamma\in\Gamma} J_\gamma(\mathbf{B}_{V_\gamma}) ) \right)} \\
&= \overline{ \amconv_\C  \left(\bigcup_{\gamma\in\Gamma} (J_\gamma)_c(\kappa \mathbf{B}_{V_\gamma}) \right) } = \overline{ \amconv_\C  \left(\bigcup_{\gamma\in\Gamma} (J_\gamma)_c( \amconv_\C \left(\kappa \mathbf{B}_{V_\gamma}) \right)\right) } \\
&= \overline{ \amconv_\C  \left(\bigcup_{\gamma\in\Gamma} (J_\gamma)_c(  \mathbf{B}_{(V_\gamma)_c} )\right) } = \overline{\mathbf{B}}_{\ell_1\big(\{(V_\gamma)_c\}_{\gamma \in\Gamma} \big)}.
\end{align*}

\end{proof}

\section{Completely Lipschitz functions}

It is now time to introduce the main object of study of this article: completely Lipschitz functions defined with respect to a suitable matrix set $\mathbf{A} = (A_n)_{n=1}^\infty$. The translation from the Banach space concept follows a similar pattern as other classes of functions --- such as completely bounded linear mappings, Schur polynomials \cite{Defant-Wiesner}, and completely holomorphic mappings  \cite{chavez2026linearizing} --- in the sense that the operator space concept depends on the amplifications of the original function.
In our case, the idea is that a function $f : A_1 \to W$ defined on the first level of a matrix set $\mathbf{A}$ is completely Lipschitz if all the $n$-th amplifications $f_n$ are uniformly Lipschitz when restricted to the corresponding level $A_n$. Let us state the precise definition:

A matrix set $\mathbf{A} = (A_n)_{n=1}^\infty$ is said to be \emph{well-founded} if $0\in A_1$ and for every $n>1$ we have $A_n \subseteq M_n(A_1)$ (that is, the entries of the elements of the higher levels all come from the first level). Notice that in this case, for any function $f : A_1 \to X$ its amplifications $f_n : A_n \to M_n(X)$ are well defined.

\begin{definition}\label{defn-completely-Lipschitz}
Let $V$ and $W$ be operator spaces, and $\mathbf{A} = (A_n)_{n=1}^\infty$ a well-founded matrix set over $V$. 
A function $f : A_1 \to W$ is said to be a \emph{completely Lipschitz function} if there exists a constant $C \ge 0$ such that for every $n\in\N$ and every $x,y \in A_n$ we have
\[
 \n{ f_n(x)-f_n(y)}_{M_n(W)} \le C \n{x-y}_{M_n(V)}.
\]
The space of all such completely Lipschitz functions taking the value $0$ at $0$ is denoted by $\Lipcb(\mathbf{A},W)$.
The infimum of the $C$'s above is denoted by $\n{f}_{\Lipcb(\mathbf{A},W)}$ or simply $\n{f}_{\Lipcb}$. In the case $W=\K$, we write $\Lipcb(\mathbf{A})$.
\end{definition}

One might argue that our terminology above is slightly imprecise, since for a function $f : A_1 \to W$ being completely Lipschitz depends on the choice of $\mathbf{A}$.
This is a conscious choice made for the sake of simplicity: when we talk about completely Lipschitz functions we will not explicitly mention the $\mathbf{A}$ unless there is a risk of confusion.
This is akin to the standard usage of saying that a function $f : X \to \R$ is Lipschitz without explicitly specifying which metric is being considered on the set $X$.

By making different choices of matrix set $\mathbf{A}$, our definition recovers all of the previously considered notions of Lipschitz functions in the operator space setting:
\begin{enumerate}[(a)]
    \item In \cite{Braga-CD}, the terminology \emph{completely Lipschitz} is used for the choice $\mathbf{A}= ( M_n(V) )_{n=1}^\infty$.
    \item In \cite{Braga-CD-Sinclair}, the name \emph{$N$-Lipschitz} is used for the choice $\mathbf{A}= ( M_n(V) )_{n=1}^N$ (or more generally, $\mathbf{A}= ( M_n(A_1) )_{n=1}^N$ for some fixed $A_1$).
    \item In \cite{braga2024small} completely Lipschitz functions with $\mathbf A=\mathbf B_V$ are called \emph{completely Lipschitz in small scale} functions.
    The choice $\mathbf{A}= ( M_n(B_V) )_{n=1}^\infty$ is also considered there.
\end{enumerate}
Furthermore, even the classical case can be recovered with our definition: if $X$ is a pointed metric space, by realizing it inside its Lipschitz-free space $\F(X)$, and letting $V$ be an operator space structure on $\F(X)$ which is compatible with the norm, the choice $A_1 = X$, $A_n = \emptyset$ for $n>1$ recovers the classical definition of a Lipschitz function.

Note that a more concise way of writing the norm of a completely Lipschitz function $f : A_1 \to W$ is
\[
\n{f}_{\Lipcb(\mathbf{A},W)} = \sup\Big\{ \n{f_{n}(x)-f_n(y)}_{M_{n}(W)}/\n{x-y}_{M_{n}(V)} \;:\; n\in\N, x,y \in A_n, x\not=y \Big\}.
\]
We will later frequently use the set of pairs appearing in the formula above, and we denote it by
$P_{\mathbf{A}} = \{ (x,y) \; : \; n\in\N, x,y \in A_n, x\not=y \}$. 
Then, the map $\Phi_{\mathbf{A},W} : \Lipcb(\mathbf{A},W) \to \ell_\infty(\{ M_{n_x}(W) \,:\, (x,y) \in P_{\mathbf{A}} \})$  given by $\Phi_{\mathbf{A},W}(f) = \big( \big(f_{n_x}(x)-f_{n_x}(y)\big)/\n{x-y} \big)_{(x,y) \in P_{\mathbf{A}}}$ is an isometry, where for each $x \in \mathbf{A}$ the number $n_x\in\N$ is the matrix level to which $x$ belongs (that is, $x \in M_{n_x}(V)$).
Note that for $(x,y) \in P_{\mathbf{A}}$, $x$ and $y$ belong to the same level of $\mathbf{A}$, and thus  $n_x = n_y$.
Moreover, since $\ell_\infty(\{ M_{n_x}(W) \,:\, (x,y) \in P_{\mathbf{A}} \})$ has a natural operator space structure, identifying $\Lipcb(\mathbf{A},W)$ with its image under $\Phi_{\mathbf{A},W}$ we can define norms for matrices of elements of $\Lipcb(\mathbf{A},W)$. Specifically, for $(f_{rs}) \in M_n( \Lipcb(\mathbf{A},W) )$ we have that
\begin{multline}\label{eqn-oss-for-Hcb}
\n{ (f_{rs}) }_{M_n(\Lipcb(\mathbf{A},W))} \\= 
\sup \left\{ \frac{ \n{\big( f_{rs}(x_{jk})-f_{rs}(y_{jk})\big)}_{M_{mn}(W)} }{\n{(x_{jk}-y_{jk})}_{M_m(V)}} \;:\; m\in\N, (x_{jk}),(y_{jk}) \in A_m, (x_{jk})\not=(y_{jk}) \right\}.
\end{multline}

Next we see that with this structure, $\Lipcb(\mathbf{A},W)$ is an operator space.
Note that $\Lipcb(\mathbf{A},W)$ being a real or complex operator space depends only on whether $W$ is real or complex, and in either case we only use the real structure on $V$ (even less than that, only the ``metric'' structure).

\begin{lemma}\label{lemma-Hcb-Banach}
For any operator spaces $V$ and $W$,
and a well-founded matrix set $\mathbf{A}= (A_n)_{n=1}^\infty$ over $V$,
the space
$\Lipcb(\mathbf{A},W)$ is an operator space.    
\end{lemma}

\begin{proof}
It is clear that $\n{\cdot}_{M_n(\Lipcb(\mathbf{A},W))}$ are indeed norms satisfying Ruan's conditions, so we only need to show completeness. Let $(f^{[n]})$ be a Cauchy sequence in $\Lipcb(\mathbf{A},W)$, so in particular it is bounded: there exists $C \ge 0$ such that for each $n\in\N$ we have $\n{f^{[n]}}_{\Lipcb} \le C$. 
Note that $(f^{[n]})$ is Cauchy in $\Lip_0(A_1,W)$, so it has a limit in there which is a Lipschitz function $f : A_1 \to W$ and in particular $(f^{[n]})$ converges to $f$ pointwise.
For each fixed $m\in\N$ and $x,y \in A_m$,
since for all $n\in\N$ we have
\[
 \n{\big( f^{[n]}_m(x)-f^{[n]}_m(y)\big)}_{M_m(W)} \le C \n{x-y}_{M_m(V)},
\]
taking the pointwise limit we conclude
\[
 \n{\big( f_m(x)-f_m(y)\big)}_{M_m(W)} \le C \n{x-y}_{M_m(V)}
\]
and therefore $f \in \Lipcb(\mathbf{A},W)$.
A standard argument shows that $(f^{[n]})$ converges to $f$ in $\Lipcb(\mathbf{A},W)$.
\end{proof}

From now on, when we say that $V$ and $W$ are operator spaces we mean that both are operator spaces over the \emph{same} field $\R$ or $\C$, unless we state it otherwise.

\begin{remark}\label{remark-linear-cb-is-holomorphic-cb}
It is clear that if $T\in\CB(V,W)$,
and $\mathbf{A}= (A_n)_{n=1}^\infty$ is a well-founded matrix set over $V$,
then the restriction of $T$ to $A_1$ is a completely Lipschitz function.
Furthermore, this restriction process yields a complete contraction $\CB(V,W) \to \Lipcb(\mathbf{A},W)$.
If we additionally assume that for every $n\in\N$ we have $\inte A_n \not=\emptyset$ (which we will shorten to $\inte \mathbf{A} \not=\emptyset$), then we in fact get a completely isometric embedding (because for a linear map between normed spaces, its restriction to any nonempty open ball has the same Lipschitz constant as the whole linear map).
\end{remark}

\begin{examples}\label{examples completely Lip}
Let us see some examples of completely Lipschitz and non completely Lipschitz functions with the main focus on the case $\mathbf A=\mathbf B_V$ (the matrix unit ball). 
\begin{enumerate}[(i)]
    \item The conjugation $f(z)=\conj( z)$ clearly belongs to $\Lipcb(\mathbf{B}_\C)$ and $\|f\|_{\Lipcb}=1$. Moreover, $f$ belongs to $\Lipcb(\mathbf{A})$ for any well-founded matrix set $\mathbf A$ over $\C$. Furthermore, given any real operator space $V$ with reasonable complexification $V_c$, the function $f(v)=\conj(v)$ belongs to $\Lipcb(\mathbf{A},V_c)$ for any well-founded matrix set $\mathbf A$ over $V_c$. 
    \item By \cite[Thm. 5]{braga2024small}, any polynomial in one complex variable without constant term belongs to $\Lipcb(\mathbf{B}_\C)$; see Proposition \ref{prop-Lipcb-is-algebra} below for a related result.
    More generally, since in a finite-dimensional space all polynomials are of finite type, Proposition \ref{prop-Lipcb-is-algebra} also shows that when $\dim(V)<\infty$ any scalar-valued polynomial on $V$ without constant term belongs to $\Lipcb(\mathbf{B}_V)$.
    \item 
    Note that \cite[Thm. 16]{braga2024small} provides a general recipe for constructing examples in $\Lipcb(\mathbf{B}_V,W)$: compositions of the form $Q \circ \restr{T}{B_V}$ where $T : V \to V_1 \oplus_\infty \cdots \oplus_\infty V_m$ is a completely bounded linear map, and $Q : V_1 \oplus_\infty \cdots \oplus_\infty V_m \to W$ is an $m$-linear map satisfying $\sup_n \n{Q_n} < \infty$, where the norm is the standard norm for an $m$-linear map (the supremum of the norms of evaluations on $m$-tuples where each element has norm at most one), and $Q_n : M_n(V_1) \oplus_\infty \cdots \oplus_\infty M_n(V_m) \to M_n(W)$ denotes the $n$-th amplification of $Q$ given by
    \[
    Q_n \big( (v^1_{jk})_{jk}, \dotsc, (v^m_{jk})_{jk} \big) = \big( Q(v^1_{jk}, \dotsc, v^m_{jk} ) \big)_{jk}.
    \]
    Such $m$-linear maps are called \emph{completely controlled} in \cite{braga2024small}, and \emph{Schur completely bounded} in \cite{Defant-Wiesner}; see \cite[Prop. 18, Cors. 19, 20, and 21]{braga2024small} for specific examples.
    In particular, from the polarization results in \cite[Sec. 7]{Defant-Wiesner} it follows that an $m$-homogeneous Schur polynomial $P:V\to W$ in the sense of \cite{Defant-Wiesner} (that is, satisfying $\sup \left\{ \n{\big(P(x_{jk})\big)}_{M_n(W)} \;:\; n\in\N, (x_{jk}) \in B_{M_n(V)} \right\} <\infty $) belongs to $\Lipcb(\mathbf{B}_V,W)$. 
   \item\label{examples completely Lip:modulo} It is shown in \cite[Prop. 3.1]{Braga-CD}  that the modulus function $|\cdot|$ does not belong to $\Lipcb(\mathbf{B}_\R)$. 
      From that, easy algebraic manipulations show that for constants $a,b\in\R$ the function $f:B_{\R}\to\R$ given by 
   \[
f(x) = \begin{cases} a x &x\ge 0\\
bx  &x<0\\
\end{cases}
\]
 is not in $\Lipcb(\mathbf{B}_\R)$ unless $a=b$.
As a consequence, for any operator space $V\not=\{0\}$ we have that $\Lipcb(\mathbf{B}_V)$ and $\Lipo(B_V)$ are different as sets. Indeed, let $v_0\in V$ be a vector of norm one and $v' \in V'$ a norming functional for $v_0$. The function $f(v) = |v'(v)|$ is clearly in $\Lipo(B_V)$, because it is a composition of Lipschitz functions, but $f$ is not in $\Lipcb(\mathbf{B}_V)$ since its restriction to $\R v_0$ reduces to the modulus.
 \item Changing the matrix set, for any operator space $V$, by \cite[Thm. 1.1]{Braga-CD} (resp. \cite[Thm. 3]{braga2024small}), the space $\Lipcb\big((M_n(V))_{n=1}^\infty, W\big)$ (resp. $\Lipcb\big((M_n(B_V))_{n=1}^\infty, W\big)$) consists exactly of the real-linear completely bounded maps $V \to W$ (resp. their restrictions to $B_V$). 
\end{enumerate}
\end{examples}

\begin{remark}
    The arguments from \cite[Prop. 3.1]{Braga-CD} also allow us to prove the following: if $f \in \Lipcb(\mathbf{B}_V,W)$ is an even function (that is, $f(-x)=f(x)$, for all $x\in B_V$), then for any $x \in B_V$
\[
\n{f(x)} \le \sqrt{2} \n{f}_{\Lipcb(\mathbf{B}_V,W)} \n{x}^2.
\]
In particular, $f$ is Fr\'echet differentiable at $0$ and its derivative is $0$.

\rm Indeed, for each $k\in\N \cup \{0\}$, let $\mathbbm{1}_{2^k}\in M_{2^k}$ be the matrix whose entries are all 1 and let $R_{2^k}\in M_{2^k}$ be the Hadamard matrix given by the following classical construction of Sylvester:
\[
R_1 = 1, \qquad \qquad R_{2^{k+1}}=\begin{bmatrix}
R_{2^k} & R_{2^k}\\
R_{2^k} & -R_{2^k}\\
\end{bmatrix}.
\] It is known that $\|R_{2^k}\|=2^{k/2}$ and $\|\mathbbm{1}_{2^k}\|= 2^k$, for all $k$.
For a fixed $x \in B_V$, choose $k$ to be the maximal nonnegative integer for which $2^{k/2}\n{x} < 1$.
Since $f$ is even, we have
\begin{equation}\label{eqn-diff-at-zero-1}
    \left\|f_{2^k}\left(R_{2^k}\otimes x\right)\right\| = \left\|f_{2^k}\left(\mathbbm{1}_{2^k}\otimes x\right)\right\|=\left\|\mathbbm{1}_{2^k}\otimes f\left(x\right) \right\|=\|\mathbbm{1}_{2^k}\|\cdot \left\|f\left(x\right) \right\|= 2^k \left\|f\left(x\right) \right\|,
\end{equation}
where $f_{2^k}$ is the $2^k$-amplification of $f$.
Now, since $\n{R_{2^k}\otimes x} = 2^{k/2}\n{x}<1$ and $f$ vanishes at $0$, the completely Lipschitz condition implies
\begin{equation}\label{eqn-diff-at-zero-2}
\left\|f_{2^k}\left(R_{2^k}\otimes x\right)\right\| \le \|f\|_{\Lipcb} \left\|R_{2^k}\otimes x\right\| = \|f\|_{\Lipcb} 2^{k/2} \|x\|.
\end{equation}
Together, \eqref{eqn-diff-at-zero-1} and \eqref{eqn-diff-at-zero-2} imply
\[
\n{f(x)} \le \|f\|_{\Lipcb} \frac{1}{2^{k/2}} \|x\| \le  \|f\|_{\Lipcb} \sqrt{2}\| x\|^2,
\]
where in the last step we have used that the maximality of $k$ implies $2^{k/2}\sqrt{2}\|x\| \ge 1$.
\end{remark}

In a subsequent article \cite{CD-Dimant-futuro} we present a precise description of $ \Lipcb(\mathbf B_{\R})$ which, in particular, shows that each function in this space (not just the even ones) has derivative at 0.

Our next result is the operator space version of the fact that the product of two bounded Lipschitz functions is again Lipschitz \cite[Prop. 1.30.(i)]{Weaver-Lipschitz-algebras}.
Closely related results can be found in \cite[Prop. 15]{braga2024small} and \cite[Prop. 3.5]{chavez2026linearizing}.

First, some notation. 
If $\mathbf{A}= (A_n)_{n=1}^\infty$ is a well-founded matrix set over an operator space $V$, and $W$ is an operator space, we let $L_\infty(\mathbf{A},W)$ be the space of functions $f : A_1 \to W$ such that $\sup \{ \n{f_n(x)}_{M_n(W)} \; :\; n\in\N, x\in A_n \}$ is finite, and denote this quantity by $\n{f}_{L_\infty(\mathbf{A},W)}$.
The identification $M_N(L_\infty(\mathbf{A},W)):= L_\infty(\mathbf{A},M_N(W))$ even allows us to endow $L_\infty(\mathbf{A},W)$ with an operator space structure, although we will not be making use of it.

\begin{proposition}\label{prop-Lipcb-is-algebra}
Suppose $\mathbf{A}= (A_n)_{n=1}^\infty$ is a well-founded matrix set over an operator space $V$, and $\mathcal{X}$ is an operator algebra.
\begin{enumerate}[(a)]
    \item If 
    If $f,g \in \Lipcb(\mathbf{A},\mathcal{X}) \cap L_\infty(\mathbf{A},\mathcal{X})$ (resp. $f,g \in L_\infty(\mathbf{A},\mathcal{X})$), then their (pointwise) product $fg$ also belongs to $\Lipcb(\mathbf{A},\mathcal{X})$ (resp. $L_\infty(\mathbf{A},\mathcal{X})$) and moreover
    \[
    \n{fg}_{\Lipcb(\mathbf{A},\mathcal{X})} \le \n{f}_{L_\infty(\mathbf{A},\mathcal{X})} \n{g}_{\Lipcb(\mathbf{A},\mathcal{X})} + \n{f}_{\Lipcb(\mathbf{A},\mathcal{X})}\n{g}_{L_\infty(\mathbf{A},\mathcal{X})},
    \]
    (resp. $\n{fg}_{L_\infty(\mathbf{A},\mathcal{X})} \le \n{f}_{L_\infty(\mathbf{A},\mathcal{X})} \n{g}_{L_\infty(\mathbf{A},\mathcal{X})}$).
    \item Suppose that additionally $\mathbf{A}$ is bounded and satisfies that for every $n\in \N$ with $A_n \not=\emptyset$ we have $0 \in A_n$.
    Then, the (pointwise) product extends to a completely bounded map $\Lipcb(\mathbf{A},\mathcal{X}) \otimes_h \Lipcb(\mathbf{A},\mathcal{X}) \to \Lipcb(\mathbf{A},\mathcal{X})$, and thus $\Lipcb(\mathbf{A},\mathcal{X})$ is completely isomorphic to an operator algebra.
\end{enumerate}
\end{proposition}

\begin{proof}
\emph{(a)} Recall a classical inequality of Schur \cite[Satz III]{Schur} (see also \cite{Ong}), which asserts that the norm of the Schur product of two matrices in $M_n$ (i.e., their entrywise product) does not exceed the product of their individual norms.
Now, for $N\in\N$ let us denote by $\square : M_N(\mathcal{B}(H)) \times M_N(\mathcal{B}(H)) \to M_N(\mathcal{B}(H))$ the (block) Schur product, that is, for any $(A_{pq}),(B_{pq}) \in M_N(\mathcal{B}(H))$ define $(A_{pq}) \square (B_{pq}) = (A_{pq}B_{pq})$.
It is shown in \cite[Thm. 2.3]{Christensen} that $\square$ is a completely bounded bilinear map of norm 1, for any complex Hilbert space $H$.
Note that the same result is immediately true for a real Hilbert space $H$, since the canonical inclusion $\mathcal{B}(H) \to \mathcal{B}(H_c)$ is a completely isometric homomorphism.

Realizing $\mathcal{X}$ as a subalgebra of $\mathcal{B}(H)$ we get that for every $(x_{pq}),(y_{pq}) \in A_N$, using the contractivity of $\square$,
\begin{align*}
    \n{\big( fg(x_{pq})-fg(y_{pq})\big)}_{M_N(\mathcal{X})} &\le   \n{\big(f(x_{pq})[g(x_{pq})-g(y_{pq})] \big)}_{M_N(\mathcal{X})} +  \n{\big([f(x_{pq})-f(y_{pq})] g(y_{pq})\big)}_{M_N(\mathcal{X})}\\
    &\le \n{\big(f(x_{pq})\big)}_{M_N(\mathcal{X})} \n{\big(g(x_{pq})-g(y_{pq}) \big)}_{M_N(\mathcal{X})} \\ &\quad+  \n{\big(f(x_{pq})-f(y_{pq})\big) }_{M_N(\mathcal{X})} \n{\big(g(y_{pq})\big)}_{M_N(\mathcal{X})} \\
    &\le  \big( \n{f}_{L_\infty(\mathbf{A},\mathcal{X})} \n{g}_{\Lipcb(\mathbf{A},\mathcal{X})} + \n{f}_{\Lipcb(\mathbf{A},\mathcal{X})}\n{g}_{L_\infty(\mathbf{A},\mathcal{X})} \big) \n{\big( x_{pq}-y_{pq}\big)}_{M_N(V)},
\end{align*}
which yields the desired conclusion.
It is clear that the same type of argument shows the ``respectively'' part.

\emph{(b)}
Since $\mathbf{A}$ is bounded, there exists a constant $C$ such that for all $N\in \N$ and all $x = (x_{pq}) \in A_N$ we have $\n{x}_{M_N(V)} \le C$. If such an $x$ exists then $0 \in A_N$, so it follows that for any $(f_{jk}) \in M_n(\Lipcb(\mathbf{A},\mathcal{X}))$
\begin{equation}
\begin{aligned}\label{eqn-prop-Lipcb-is-algebra}
    \n{\big(f_{jk}(x_{pq})\big)}_{M_{nN}(\mathcal{X})} &= \n{\big(f_{jk}(x_{pq}) - f_{jk}(0)\big)}_{M_{nN}(\mathcal{X})} \\
    &\le \n{(f_{jk})}_{M_n(\Lipcb(\mathbf{A},\mathcal{X}))} \n{(x_{pq})}_{M_N(V)} \le C \n{(f_{jk})}_{M_n(\Lipcb(\mathbf{A},\mathcal{X}))}.
\end{aligned}
\end{equation}

Therefore, it follows from the previous part that the pointwise product defines a bilinear map
\[
\Lipcb(\mathbf{A},\mathcal{X}) \times \Lipcb(\mathbf{A},\mathcal{X}) \to \Lipcb(\mathbf{A},\mathcal{X}).
\]
Let $f = (f_{jk}) \in M_{m,n}(\Lipcb(\mathbf{A},\mathcal{X}))$ and 
$g = (g_{k \ell}) \in M_{n,m}(\Lipcb(\mathbf{A},\mathcal{X}))$.
By the description of the Haagerup tensor product as in \cite[Eqn. (9.2.1)]{Effros-Ruan-book}, we need to estimate the norm of $\big(\sum_{k=1}^n f_{jk}g_{k\ell}\big)_{j\ell}$ in $M_{m}(\Lipcb(\mathbf{A},\mathcal{X}))$.
To that effect, consider $N\in\N$ and $(x_{pq}),(y_{pq}) \in A_N$. We then have
\begin{equation}
\begin{aligned}\label{eqn-Haagerup-product-1}
     \n{\big( \sum_{k=1}^n f_{jk}(x_{pq})g_{k \ell}(x_{pq})-f_{jk}(y_{pq})g_{k \ell}(y_{pq}) \big)}_{M_{mN}(\mathcal{X})} 
    &\le  
    \n{\big( \sum_{k=1}^n f_{jk}(x_{pq})[g_{k\ell}(x_{pq})-g_{k\ell}(y_{pq})] \big)}_{M_{mN}(\mathcal{X})} \\
    &\quad +  \n{\big( \sum_{k=1}^n [f_{jk}(x_{pq})-f_{jk}(y_{pq})] g_{k\ell}(y_{pq})\big)}_{M_{mN}(\mathcal{X})}.
\end{aligned}    
\end{equation}

and therefore, once again using that $\square$ is a completely bounded bilinear map of norm 1,
\begin{multline}\label{eqn-Haagerup-product-2}
\n{\Big( \sum_{k=1}^n \big(f_{jk}(x_{pq})\big)_{pq} \square \big( g_{k\ell}(x_{pq})-g_{k\ell}(y_{pq})  \big)_{pq} \Big)_{j\ell}}_{M_{m}(M_N(\mathcal{X}))} \\
\le 
\n{\big(f_{jk}(x_{pq})\big)}_{M_{mN,nN}(\mathcal{X})} \n{\big(g_{k\ell}(x_{pq})-g_{k\ell}(y_{pq}) \big)}_{M_{nN,mN}(\mathcal{X})}
\end{multline} 
and analogously
\begin{multline}\label{eqn-Haagerup-product-3}
\n{\Big( \sum_{k=1}^n \big(f_{jk}(x_{pq})-f_{jk}(y_{pq})\big)_{pq} \square \big( g_{k\ell}(y_{pq})\big)_{pq} \Big)_{j\ell}}_{M_{m}(M_N(\mathcal{X}))} \\
\le 
\n{\big(f_{jk}(x_{pq})-f_{jk}(y_{pq})\big)}_{M_{mN,nN}(\mathcal{X})} \n{\big(g_{k\ell}(y_{pq}) \big)}_{M_{nN,mN}(\mathcal{X})}.
\end{multline} 

From \eqref{eqn-prop-Lipcb-is-algebra}, \eqref{eqn-Haagerup-product-1}, \eqref{eqn-Haagerup-product-2} and \eqref{eqn-Haagerup-product-3} we have
\[
\n{\big(\sum_{k=1}^n f_{jk}g_{k\ell}\big)_{j\ell}}_{M_{m}(\Lipcb(\mathbf{A},\mathcal{X}))} \le 2C \n{(f_{jk})}_{M_{m,n}(\Lipcb(\mathbf{A},\mathcal{X}))} \n{(g_{k\ell})}_{M_{n,m}(\Lipcb(\mathbf{A},\mathcal{X}))},
\]
from where we conclude that the (pointwise) product extends to a completely bounded map $\Lipcb(\mathbf{A},\mathcal{X}) \otimes_h \Lipcb(\mathbf{A},\mathcal{X}) \to \Lipcb(\mathbf{A},\mathcal{X})$ of $\cb$-norm at most $2C$.

By \cite[Thm. 5.2.1]{BlecherLeMerdy2004} (note that \cite[p. 31]{blecher2023real} remarks that the real case follows from the complex one by complexification), if $\mathcal{A}$ is an algebra which is also an operator space, then the multiplication extends to  a completely bounded map $\mathcal{A} \otimes_h \mathcal{A} \to \mathcal{A}$  if and only if $\mathcal{A}$ is completely isomorphic to an operator algebra. This yields the last part of the desired conclusion.
\end{proof}

The reader might have noticed that our proof of part $(b)$ in the previous result can be used to obtain a generalization of part $(a)$. We have chosen clarity over generality, to emphasize how part $(a)$ generalizes the classical result \cite[Prop. 1.30.(i)]{Weaver-Lipschitz-algebras}.

As happens in the Banach (and complete metric) case, Lipschitz functions on a set extend to the completion with the same Lipschitz norm \cite[Prop. 1.6]{Weaver-Lipschitz-algebras}.

\begin{lemma}\label{lemma:extension al borde}
    Let $V$ and $W$ be operator spaces, $\mathbf{A}= (A_n)_{n=1}^\infty$ a well-founded matrix set over $V$,  and $f\in\Lipcb(\mathbf{A},W)$. Then, $f$ can be extended to $\overline{A_1}$ and satisfies that for every $n\in\N$ and every $(x_{jk}),(y_{jk}) \in  \overline{A_n}$,
\[
 \n{\big( f(x_{jk})-f(y_{jk})\big)}_{M_n(W)} \le \|f\|_{\Lipcb}\n{\big( x_{jk}-y_{jk}\big)}_{M_n(V)},
\]
that is, the extension belongs to $\Lipcb(\overline{\mathbf{A}},W)$.
\end{lemma}

\begin{proof}
   The result follows just applying \cite[Prop. 1.6]{Weaver-Lipschitz-algebras} to each matrix level. 
\end{proof}

In the classical case, for a metric space $X$ it is known that the space $\Lipo(X,\C)$ can be identified real-isometrically with the Taylor complexification of $\Lipo(X,\R)$ \cite[Prop. 2.4]{Abbar-Coine-Petitjean}. Our next result is a far more general operator space version of this.

\begin{proposition}\label{prop-Lipcb-is-complexification}
Let $V$ be an operator space and let $W$ be a real operator space.
Let $\mathbf{A}= (A_n)_{n=1}^\infty$ be a well-founded matrix set over $V$.
Then the complex operator space $\Lipcb(\mathbf{A},W_c)$ is the unique reasonable operator space complexification of the real operator space $\Lipcb(\mathbf{A},W)$.
\end{proposition}

\begin{proof}
If $f : A_1 \to W$ is completely Lipschitz, it is clear that $\kappa_W \circ f : A_1 \to W_c$ is also completely Lipschitz, since $\kappa_W$ is a (real) complete isometry.
In fact, the assignment $f \mapsto \kappa_W \circ f$ yields a real complete isometry $\kappa : \Lipcb(\mathbf{A},W) \to \Lipcb(\mathbf{A},W_c)$.

Note also that $f : A_1 \to W_c$ is completely Lipschitz if and only if $\Re(f), \Im(f) : A_1 \to W$ are  completely Lipschitz:
the ``only if'' part follows from the fact that $\Re, \Im : W_c \to W$ are real complete contractions (Corollary \ref{cor-real-part-is-complete-contraction}), whereas
the ``if'' is clear from $f = \Re(f) + i\Im(f)$.
Therefore,
\[
\Lipcb(\mathbf{A},W_c) = \Lipcb(\mathbf{A},W) + i \Lipcb(\mathbf{A},W)
\]
with the first summand being equal to the image of the map $\kappa$ above; this shows that $\Lipcb(\mathbf{A},W_c)$ is a vector space complexification of $\Lipcb(\mathbf{A},W)$.
Finally, note that for any $(f_{jk}),(g_{jk}) \in M_n(\Lipcb(\mathbf{A},W))$ we have the reasonable condition
\[
\n{(f_{jk} + ig_{jk})_{jk}}_{M_n(\Lipcb(\mathbf{A},W_c))} = \n{(f_{jk} - ig_{jk})_{jk}}_{M_n(\Lipcb(\mathbf{A},W_c))}, 
\]
which is inherited from the reasonable condition for $W_c$, and thus it follows that $\Lipcb(\mathbf{A},W_c)$ is the unique reasonable operator space complexification of $\Lipcb(\mathbf{A},W)$ \cite[Thm. 3.1]{Ruan-complexifications}.
\end{proof}

\section{The predual of $\Lipcb(\mathbf{A})$: Dixmier--Ng--Kaijser approach}
\label{sec-Dixmier-Ng-Kaijser}

A very convenient way to verify that a Banach space is a dual space is by using the Dixmier--Ng theorem, which has the additional advantage of explicitly identifying a predual. In the operator space setting, Le Merdy \cite{LeMerdy-duality} showed that having a predual as a Banach space does not imply having one as an operator space, and provided necessary and sufficient conditions for a positive answer.
The complex version of the following proposition appeared in \cite[Prop. 3.1]{LeMerdy-duality}.

\begin{proposition}\label{prop-LeMerdy}
Let $V$ be a real operator space, and let $Y$ be a Banach space such that there exists an isometric identification $\iota : V \to Y'$ as Banach spaces. Endow $Y$ with its predual real operator space structure, that is, the one coming from the inclusion $Y \hookrightarrow Y''=V'$. The following are equivalent:
\begin{enumerate}[(a)]
    \item $\iota : V \to Y'$ is a complete isometry.
    \item For each $n\in\N$, $\overline{B}_{M_n(V)}$ is closed for the entrywise $\sigma(V,Y)$-topology. 
    \item For each $n\in\N$, $\overline{B}_{M_n(Y)}$ is dense in $\overline{B}_{M_n(V')}$ for the entrywise $\sigma(V',V)$-topology. 
    \item For each $n\in\N$, $M_n(Y)'' = M_n(V')$ isometrically. 
\end{enumerate}
\end{proposition}

\begin{proof}

The proof of the complex case in \cite[Prop. 3.1]{LeMerdy-duality} works verbatim in the real case. 
The basic properties of the Haagerup tensor product needed in the argument also hold in the real case, see \cite[Lemma 5.4]{blecher2023real} (and the discussions before and after it).

\end{proof}

A version of Kaijser \cite{Kaijser} of the Dixmier--Ng theorem  combined with  Le Merdy's result  about conditions for operator space preduals, yields the next theorem.
The complex version appeared in \cite[Thm. 4.1]{chavez2026linearizing}, and the same proof works in the real case using Proposition \ref{prop-LeMerdy}.

\begin{theorem}\label{thm-Kaijser}
Let $V$ be an operator space, and let $Y \subset V'$ be a set of continuous linear functionals on $V$ such that:
\begin{enumerate}[(a)]
\item $Y$ separates points of $V$.
\item For each $n\in\N$, the closed unit ball of $M_n(V)$ is compact for the entrywise $\sigma(V,Y)$-topology.
\end{enumerate}
Then $V$ is a dual operator space, and an operator space predual for $V$ is the closure in $V'$ of the linear span of the set $Y$.    
\end{theorem}

Before stating the duality theorem, we remark that from the definition of the norm on $\Lipcb(\mathbf{A})$ it is clear that for any $x\in A_1$ the evaluation functional $\delta_{\mathbf{A}}(x)$ (given by $\delta_{\mathbf{A}}(x)(f) = f(x)$ for $f\in \Lipcb(\mathbf{A})$)  belongs to $\Lipcb(\mathbf{A})'$ and   $\|\delta_{\mathbf{A}}(x)\| \le \|x\|$, since for any $f\in\Lipcb(\mathbf{A})$ we have
$|\delta_{\mathbf{A}}(x)(f)| = |f(x)| = |f(x)-f(0)| \le \n{f}_{\Lipcb} \n{x-0}$. 
Moreover, 
$\|\delta_{\mathbf{A}}(x)\|\ge |x'(x)|$ for all $x'\in \overline B_{V'}$ by Remark \ref{remark-linear-cb-is-holomorphic-cb}, hence $\|\delta_{\mathbf{A}}(x)\| = \|x\|$ for all $x\in A_1$.

\begin{theorem}\label{thm-predual-version-Dixmier-Ng-Kaijser}
Let $V$ be an operator space and
$\mathbf{A}= (A_n)_{n=1}^\infty$ a well-founded matrix set over $V$.
Then the space $\Lipcb(\mathbf{A})$ is a dual operator space, and in fact $\Fcb(\mathbf{A}) = \overline{\spa}\{ \delta_{\mathbf{A}}(x) \}_{x \in A_1} \subseteq \Lipcb(\mathbf{A})'$ is an operator space predual of $\Lipcb(\mathbf{A})$.
Moreover, on bounded subsets of $\Lipcb(\mathbf{A})$ the restriction of the weak$^*$ topology is the topology of pointwise convergence.
\end{theorem}

\begin{proof}
Let $Y =\{ \delta_{\mathbf{A}}(x) \}_{x \in A_1} \subseteq  \Lipcb(\mathbf{A})'$. It is clear that it separates points of $\Lipcb(\mathbf{A})$, and the $\sigma(\Lipcb(\mathbf{A}), Y)$-topology is the topology of pointwise convergence.
Let $(f^{[\alpha]})_{\alpha\in\Lambda}$ be a net in the closed unit ball of $\Lipcb(\mathbf{A})$. 

Note that in particular the functions $f^{[\alpha]}$ are pointwise bounded,
since for each $x \in A_1$ we have $|f(x)| = |f(x)-f(0)| \le \n{x}$.
Applying Tychonoff's theorem to $\prod_{x\in A_1} \n{x} \overline{B}_{\K}$, there is a subnet $(f^{[\gamma]})_{\gamma\in\Gamma}$ that converges pointwise to a function $f$.
For any $\gamma\in\Gamma$, any $n\in\N$ and any $x,y\in A_n$ we have $\n{f_{n}^{[\gamma]}(x)- f_n^{[\gamma]}(y)}_{M_{n}} \le \n{x-y}_{M_n(V)}$,
so it follows from the pointwise convergence that $f$ belongs to the closed unit ball of $\Lipcb(\mathbf{A})$.
Thus, $\overline B_{\Lipcb(\mathbf{A})}$ is compact in the $\sigma(\Lipcb(\mathbf{A}), Y)$-topology.

A similar argument shows that for any $n\in\N$ the closed unit ball of $M_n(\Lipcb(\mathbf{A}))$ is compact for the entrywise $\sigma(\Lipcb(\mathbf{A}), Y)$-topology, so by the Dixmier--Ng--Kaijser-Le Merdy theorem we conclude that $\Fcb(\mathbf{A})$ is an operator space predual of $\Lipcb(\mathbf{A})$.

Since $\overline B_{\Lipcb(\mathbf{A})}$ is compact in the weak$^*$ topology, the topology of pointwise convergence is Hausdorff, and the former is finer than the latter, by a standard basic result both topologies coincide on $B_{\Lipcb(\mathbf{A})}$. 
\end{proof}

When needed, we will write $\Fcb_\R(\mathbf{A})$ and $\Fcb_\C(\mathbf{A})$ to differentiate between the real and complex cases.
Now that we have found a predual for $\Lipcb(\mathbf{A})$ we show that, as in the typical linearization procedures for Banach spaces, the evaluation map $\delta_{\mathbf{A}} : A_1 \to \Fcb(\mathbf{A})$ belongs to the class of functions being linearized: it is a completely Lipschitz function (and it has norm one).

\begin{proposition}\label{prop-delta-is-cb-holomorphic}
Let $V$ be an operator space and
$\mathbf{A}= (A_n)_{n=1}^\infty$ a well-founded matrix set over $V$.
The evaluation map $\delta_{\mathbf{A}} : A_1 \to \Fcb(\mathbf{A})$ given by $x \mapsto \delta_{\mathbf{A}}(x)$
is a completely Lipschitz function, and furthermore the norm of $\delta_{\mathbf{A}}$ in $\Lipcb(\mathbf{A}, \Fcb(\mathbf{A}))$ is 1.
More precisely, $\delta_{\mathbf{A}}$ is completely isometric in the sense that for all $x,y\in A_n$ we have
\[
\n{(\delta_{\mathbf{A}})_n(x) - (\delta_{\mathbf{A}})_n(y)}_{M_n(\Fcb(\mathbf{A}))} = \n{x-y}_{M_n(V)}.
\]
\end{proposition}

\begin{proof}
Consider $n\in\N$ and $(x_{jk}),(y_{jk}) \in A_n$.
Then, from the duality $\Fcb(\mathbf{A})' = \Lipcb(\mathbf{A})$,
\begin{align*}
\n{(\delta_{\mathbf{A}}(x_{jk})-\delta_{\mathbf{A}}(y_{jk} )) }_{M_n(\Fcb(\mathbf{A}))} &= 
\sup \left\{ \n{( \delta_{\mathbf{A}}(x_{jk})f_{rs}-\delta_{\mathbf{A}}(y_{jk})f_{rs} )}_{M_{nm}} : (f_{rs}) \in B_{M_m(\Lipcb(\mathbf{A}))}, m\in\N \right\}
 \\
&=\sup \left\{ \n{( f_{rs}(x_{jk}) - f_{rs}(y_{jk}) )}_{M_{nm}} : (f_{rs}) \in B_{M_m(\Lipcb(\mathbf{A}))}, m\in\N \right\}.
\end{align*}
Therefore, we conclude from \eqref{eqn-oss-for-Hcb} that
\[
\n{(\delta_{\mathbf{A}}(x_{jk})-\delta_{\mathbf{A}}(y_{jk}) )}_{M_n(\Fcb(\mathbf{A}))} \le \n{(x_{jk}-y_{jk})}_{M_n(V)}.
\]
The reverse inequality is clear from the fact that
\[
\n{(x_{jk} - y_{jk})}_{M_n(V)} = \sup \{ \n{\big( z_{rs}'(x_{jk} - y_{jk})\big)}_{M_{nm}} \;:\; (z_{rs}') \in B_{M_m(V')} \}
\]
and Remark \ref{remark-linear-cb-is-holomorphic-cb}.
\end{proof}

Note that by Lemma \ref{lemma:extension al borde}, $\delta_{\mathbf{A}}$ extends to $\overline{A_1}$ and satisfies 
$\n{(\delta_{\mathbf{A}}(x_{jk}) - \delta_{\mathbf{A}}(y_{jk}))}_{M_n(\Fcb(\mathbf{A}))} = \n{(x_{jk}-y_{jk})}_{M_n(V)}$ for all $(x_{jk}),(y_{jk})\in \overline{A_n}$.

A classical useful property of the Lipschitz-free space is that its unit ball is the closed convex hull of the elementary molecules, see e.g. \cite[Lemma 3.2]{Abbar-Coine-Petitjean}.
In the completely Lipschitz setting, an analogous statement holds for the matrix unit ball.

\begin{proposition}\label{proposition-ball-is-convex-hull-of-image-of-delta}
Let $V$ be an operator space and
$\mathbf{A}= (A_n)_{n=1}^\infty$ a well-founded matrix set over $V$.
Let $\mathbf{E} = (E_n)_{n=1}^\infty$ be the matrix set over $\Fcb(\mathbf{A})$ given by, for each $n\in\N$,
\begin{equation}\label{eqn-image-of-delta}
    E_n = \bigg\{ \frac{(\delta_{\mathbf{A}})_n(x)-(\delta_{\mathbf{A}})_n(y)}{\n{x-y}}  \;:\; x,y \in A_n, x\not=y \bigg\}.
\end{equation}
Then the closed matrix unit ball of $\Fcb(\mathbf{A})$ is the closure of the absolutely matrix convex hull of $\mathbf{E}$.
\end{proposition}

\begin{proof}
From the definition of the operator space structure on $\Lipcb(\mathbf{A})$,
specifically \eqref{eqn-oss-for-Hcb},
note that for any $f \in M_m(\Lipcb(\mathbf{A}))$,
\begin{equation}\label{eqn-completely-norming}
\n{f}_{M_m(\Lipcb(\mathbf{A}))} = \sup \big\{ \n{\mpair{f}{e}}_{M_{nm}} \;:\; e \in E_n, n\in\mathbb N \big\}.
\end{equation}
Note that $\mathbf{E} \subseteq \overline{\mathbf{B}}_{\Fcb(\mathbf{A})}$, since by Proposition \ref{prop-delta-is-cb-holomorphic} all the elements of $E_n$ have norm one.
Therefore, the closure of the absolutely matrix convex hull of $\mathbf{E}$ is contained in $\overline{\mathbf{B}}_{\Fcb(\mathbf{A})}$. If they were different, the Hahn--Banach theorem for matrix convexity (Theorem \ref{thm-Hahn-Banach-matrix-convexity}) would yield a contradiction with \eqref{eqn-completely-norming}.
\end{proof}

The above representation yields that for an \emph{$N$-finite} well-founded matrix set $\mathbf A$ (that is, where the levels above the $N$-th one are empty) the space $\Fcb(\mathbf{A})$ is $N$-maximal.

\begin{proposition}
Let $V$ be an operator space, 
and suppose that the well-founded matrix set $\mathbf{A}= (A_n)_{n=1}^\infty$ over $V$ satisfies that $A_n = \emptyset$ for every $n>N$. Then the operator space $\Fcb(\mathbf{A})$ is $N$-maximal.
\end{proposition}

\begin{proof}
It follows from Proposition \ref{proposition-ball-is-convex-hull-of-image-of-delta} that
$\overline{\mathbf{B}}_{\Fcb(\mathbf{A})} = \overline{\amconv(\mathbf{E})}$ where the matrix set $\mathbf{E}$ is given by \eqref{eqn-image-of-delta} and therefore satisfies $E_n=\emptyset$ for all $n>N$.
Lemma \ref{lemma-convex-hull-of-finite-height} thus yields that $\Fcb(\mathbf{A})$ is $N$-maximal.

\end{proof}

The following is a predual version of Proposition \ref{prop-Lipcb-is-complexification}, relating the Lipschitz-free spaces in the real and complex cases.
It is an operator space version of \cite[Prop. 2.4]{Abbar-Coine-Petitjean}.

\begin{proposition}
Let $V$ be an operator space, 
and let $\mathbf{A}= (A_n)_{n=1}^\infty$ be a well-founded matrix set over $V$.
Then $\Fcb_\C(\mathbf{A})$ is the unique reasonable operator space complexification of $\Fcb_\R(\mathbf{A})$.
\end{proposition}

\begin{proof}
We already know from Proposition \ref{prop-Lipcb-is-complexification} that $\Lipcb(\mathbf{A},\C)$ is the unique operator space complexification of $\Lipcb(\mathbf{A},\R)$, which in turn implies that $\Lipcb(\mathbf{A},\C)'$ is the unique operator space complexification of $\Lipcb(\mathbf{A},\R)'$ (because the dual of the complexification is the complexification of the dual \cite{Ruan-complexifications}).
The inclusion $\Fcb_\R(\mathbf{A}) \hookrightarrow \Lipcb(\mathbf{A},\R)'$ is completely isometric, so by complexification \cite[Prop. 1.1]{blecher2023real} we also get a completely isometric inclusion
\begin{equation}\label{eqn-Fcb-real-and-complex}
    \Fcb_\R(\mathbf{A})_c \hookrightarrow \big(\Lipcb(\mathbf{A},\R)'\big)_c =  \Lipcb(\mathbf{A},\C)'.
\end{equation}
For an $x\in A_1$, we denote $\delta^\R_{\mathbf{A}}(x) : \Lipcb(\mathbf{A},\R) \to \R$ and $\delta^\C_{\mathbf{A}}(x) : \Lipcb(\mathbf{A},\C) \to \C$ the corresponding evaluation functionals. 
Since $\{ \delta^\R_{\mathbf{A}}(x) \}_{x \in A_1} $ $\R$-spans a dense subspace of $\Fcb_\R(\mathbf{A})$, its image under the isometric inclusion in \eqref{eqn-Fcb-real-and-complex} $\C$-spans a dense subspace of the copy of $\Fcb_\R(\mathbf{A})_c$ inside $\Lipcb(\mathbf{A},\C)'$.
It is easy to see that each $\delta^\R_{\mathbf{A}}(x)$ gets mapped to the corresponding $\delta^\C_{\mathbf{A}}(x)$, finishing the proof.
\end{proof}

\begin{remark}
In the Banach setting, an important result of Kaufmann \cite[Cor. 3.3]{Kaufmann} shows that in the real case the Lipschitz-free spaces of a Banach space $X$ and of its unit ball $B_X$ are isomorphic.
Thus, one may wonder if the analogous result holds in the operator space case. But because of the freedom we have in choosing matrix sets, one can easily force the answer to be trivial either way:
\begin{enumerate}[(a)]
    \item For a real operator space $V$, it follows from \cite[Thm. 1.1]{Braga-CD} and \cite[Thm. 3]{braga2024small} that the spaces $\Fcb\big((M_n(V))_{n=1}^\infty\big)$ and $\Fcb\big((M_n(B_V))_{n=1}^\infty\big)$ are completely isometric (and in fact completely isometric to $V$ itself).
    \item However, the spaces $\Fcb\big((M_n(V))_{n=1}^\infty\big)$ and $\Fcb\big((B_{M_n(V)})_{n=1}^\infty\big)$ are generally different: if $V\not=\{0\}$ is finite-dimensional, so is the former (since it is completely isometric to $V$), but the latter is infinite-dimensional (this follows from Proposition \ref{prop-Lipcb-is-algebra}). 
\end{enumerate}
\end{remark}

\section{The predual of $\Lipcb(\mathbf{A})$: de Leeuw's map approach}\label{sec-de-Leeuw}

Recall that by definition of the space $\Lipcb(\mathbf{A})$, the map $\Phi_{\mathbf{A},\K} : \Lipcb(\mathbf{A}) \to \ell_\infty(\{ M_{n_x} \,:\, (x,y) \in P_{\mathbf A} \})$ is a completely isometric embedding, and
the latter space is the dual of $\ell_1(\{ S_1^{n_x} \,:\, (x,y) \in  P_{\mathbf A} \})$.

\begin{lemma}\label{prop-de-Leeuws-map}
For any $\K$-operator space $V$,
and any well-founded matrix set $\mathbf{A}= (A_n)_{n=1}^\infty$ over $V$,
the image of $\Phi_{\mathbf{A},\K}$ is a weak$^*$ closed subspace of  $\ell_\infty(\{ M_{n_x} \,:\, (x,y) \in  P_{\mathbf A} \})$.
\end{lemma}

\begin{proof}
By the Kre\u{\i}n--\v{S}mulian theorem, it suffices to check that the unit ball of the image of $\Phi_{\mathbf{A},\K}$ is weak$^*$ closed.
Let $(f^{[\gamma]})_{\alpha\in\Gamma}$ be a net in the unit ball of $\Lipcb(\mathbf{A})$ so that $(\Phi_{\mathbf{A},\K}(f^{[\gamma]}))_{\gamma\in\Gamma}$ converges in the weak$^*$ topology of $\ell_\infty(\{ M_{n_x} \,:\, (x,y) \in  P_{\mathbf A}\})$ to some limit $g$.
In particular, for any $x \in A_1\setminus\{0\}$ we have that the net $\big( (f^{[\gamma]}(x)-f^{[\gamma]}(0))/\n{x-0} \big)_{\gamma\in\Gamma}$ converges, and thus
the net $(f^{[\gamma]})_{\gamma\in\Gamma}$ converges pointwise to a function $f : A_1 \to \K$. 
Now, fix $n\in\N$ and $x,y \in A_n$. Since $\big( f^{[\alpha]}_{n}(x)-f^{[\alpha]}_{n}(y)\big)_{\alpha \in\Lambda}$ is a net in $\n{x-y}_{M_n(V)}\overline{B}_{M_{n}}$ converging pointwise to $f_{n}(x)-f_n(y)$, we conclude that $\n{f_{n}(x)-f_n(y)} \le \n{x-y}_{M_n(V)}$, so $f$ belongs to the unit ball of $\Lipcb(\mathbf{A})$. Clearly $g = \Phi_{\mathbf{A},\K}(f)$, yielding the desired result.
\end{proof}

The previous lemma provides an alternative way of finding a predual for $\Lipcb(\mathbf{A})$.

\begin{corollary}\label{corollary-predual}
The space $\Lipcb(\mathbf{A})$ is a dual operator space.
Moreover, on bounded subsets of $\Lipcb(\mathbf{A})$ the weak$^*$ topology coincides with the topology of pointwise convergence.    
\end{corollary}

\begin{proof}
By the previous lemma and either \cite[Prop. 4.2.2]{Effros-Ruan-book} (complex case) or \cite[Lemma 5.2]{blecher2023real} (real case), $\Lipcb(\mathbf{A})$ is a dual operator space. Moreover, the weak$^*$ topology on  $\Lipcb(\mathbf{A})$ is the restriction (to $\Phi_{\mathbf{A},\K}(\Lipcb(\mathbf{A}))$) of the weak$^*$ topology from $\ell_\infty(\{ M_{n_x} \,:\, (x,y) \in  P_{\mathbf A}\})$.
In particular, convergence in the weak$^*$ topology of  $\Lipcb(\mathbf{A})$ implies pointwise convergence.
But on the unit ball of $\Lipcb(\mathbf{A})$ the weak$^*$ topology is compact and the topology of pointwise convergence is Hausdorff, so by standard basic topology the topologies must coincide.
\end{proof}

Note that the map $\Phi_{\mathbf{A},\K}$ is playing a role analogous to de Leeuw's map in the context of Lipschitz free spaces (see \cite[Sec. 2.4]{Weaver-Lipschitz-algebras}), hence the name of this section.

We will now show that the predual appearing in the proof above coincides with the one we had previously found with the Dixmier--Ng--Kaijser machinery.
Let $H_{\perp}$ be the preannihilator of $\Phi_{\mathbf{A},\K}(\Lipcb(\mathbf{A}))$ in $\ell_1(\{ S_1^{n_x} \,:\, (x,y) \in  P_{\mathbf A}\})$, that is,
\[
H_{\perp} = \big\{ z \in \ell_1(\{ S_1^{n_x} \,:\, (x,y) \in  P_{\mathbf A}\}) \;:\; \pair{\Phi_{\mathbf{A},\K}(f)}{z} = 0 \text{ for all } f \in \Lipcb(\mathbf{A}) \big\}.
\]
Then the predual $\mathscr{F}$ appearing in the proof of Corollary \ref{corollary-predual} is given as
\[
\mathscr{F} = \ell_1(\{ S_1^{n_x} \,:\, (x,y) \in  P_{\mathbf A}\})/H_{\perp}.
\]
Let us denote by $Q$ the associated complete metric surjection 
\[
Q :  \ell_1(\{ S_1^{n_x} \,:\, (x,y) \in  P_{\mathbf A}\}) \to \mathscr{F}.
\]
Since $\mathscr{F}' = \Lipcb(\mathbf{A})$, we can canonically completely isometrically isomorphically identify $\mathscr{F}$ with a subspace of $\Lipcb(\mathbf{A})'$. We now identify this subspace.

\begin{proposition}\label{prop-Lip-free-de-Leeuw}
For any operator space $V$, 
and any well-founded matrix set $\mathbf{A}= (A_n)_{n=1}^\infty$ over $V$,
the space $\mathscr{F}$ can be completely isometrically isomorphically identified with $\Fcb(\mathbf{A})=\overline{\spa}\{ \delta_{\mathbf{A}}(x) \}_{x \in A_1} \subseteq \Lipcb(\mathbf{A})'$.
\end{proposition}

\begin{proof}
First we will show that for each $x_0\in A_1\setminus \{0\}$, the evaluation map $\delta_{\mathbf{A}}(x_0)$ belongs to the subspace of $\Lipcb(\mathbf{A})'$ which is identified with $\mathscr{F}$.
To that end, let $e_{x_0} \in \ell_1(\{ S_1^{n_x} \,:\, (x,y) \in  P_{\mathbf A}\})$ be the element that has a $\n{x_0}$ in the $(x_0,0)$ coordinate (note that the space $S_1^{n_{x_0}}$ is one-dimensional) and $0$'s everywhere else. Then for every $f \in \Lipcb(\mathbf{A})$,
\[
\pair{Qe_{x_0}}{f} = \pair{e_{x_0}}{\Phi_{\mathbf{A},\K}(f)} = \n{x_0} \frac{f(x_0)-f(0)}{\n{x_0}} = f(x_0),
\]
that is, $Qe_{x_0}$ corresponds to $\delta_{\mathbf{A}}(x_0)$.
Therefore, additionally observing that $\delta_{\mathbf{A}}(0)=0$, we conclude that $\spa\{ \delta_{\mathbf{A}}(x) \}_{x \in A_1}$ is contained in the subspace of $\Lipcb(\mathbf{A})'$ which is identified with $\mathscr{F}$.

Next, observe that the finitely supported vectors are dense in $\ell_1(\{ S_1^{n_x} \,:\, (x,y) \in  P_{\mathbf A}\})$ and therefore their images under $Q$ are dense in $\mathscr{F}$.
Thus, it suffices to check that if we fix $n\in\N$ and $x_0=(x_{jk})_{j,k=1}^{n},y_0=(y_{jk})_{j,k=1}^{n}$ in  $A_n$ with $x_0 \not= y_0$,
and a matrix $ c = (c_{jk})_{j,k=1}^{n} \in S_1^{n}$, the element $z \in \ell_1(\{ S_1^{n_x} \,:\, (x,y) \in  P_{\mathbf A}\}) $ having $c$ in the $(x_0,y_0)$-coordinate and $0$'s everywhere else gets mapped by $Q$ inside $\mathscr{F}$. Now, for every $f \in \Lipcb(\mathbf{A})$ we have
\[
\pair{Qz}{f} = \pair{z}{\Phi_{\mathbf{A},\K}(f)} = \frac{1}{\n{x_0-y_0}_{M_n(V)}}\sum_{j,k=1}^{n} c_{jk} \big( f(x_{jk})- f(y_{jk}) \big)
\]
that is, the element of $\Lipcb(\mathbf{A})'$ corresponding to $Qz$ is $ \sum_{j,k=1}^{n} c_{jk} 
\big( \delta_{\mathbf{A}}(x_{jk}) - \delta_{\mathbf{A}}(y_{jk})\big)$, which belongs to $\Fcb(\mathbf{A})$. 
The desired conclusion follows.
\end{proof}

\section{The predual of $\Lipcb(\mathbf{A})$: the Arens--Eells approach}\label{sec-Arens-Eells}

The Dixmier--Ng--Kaijser approach to $\Fcb(\mathbf{A})$ from Section \ref{sec-Dixmier-Ng-Kaijser}, which is very similar to how the Lipschitz-free space is introduced in \cite{Godefroy-Kalton}, identifies $\Fcb(\mathbf{A})$ with a subspace of a dual space (namely $\Lipcb(\mathbf{A})'$).
Therefore, the norm of a matrix in $M_n(\Fcb(\mathbf{A}))$ is naturally calculated as a supremum.
Section \ref{sec-de-Leeuw} gives an alternative expression where $\Fcb(\mathbf{A})$ is identified with a quotient, and therefore the norm of a matrix in $M_n(\Fcb(\mathbf{A}))$ is now expressed as an infimum, but the expressions involved are rather unwieldy.
The purpose of this section is to present a third approach reminiscent of how the Lipschitz-free space is introduced in \cite{Weaver-Lipschitz-algebras} and which goes back to \cite{Arens-Eells},
according to which the norm of an element of the Lipschitz-free space over a metric space $X$ which is a linear combination of evaluation functionals can be calculated as
\begin{equation}\label{eqn-classical-AE-incomplete}
\n{u}_{\F(X)} = \inf \Big\{ \sum_{j=1}^n |c_j| d(x_j,y_j) \;:\; u = \sum_{j=1}^n c_j\big( \delta_X(x_j) - \delta_X(y_j) \big) \Big\}.
\end{equation}
Our next result is a noncommutative version of this.
It is closely related to Proposition \ref{proposition-ball-is-convex-hull-of-image-of-delta} above.

\begin{proposition}\label{prop-Arens-Eells-incomplete-without-matrix}
Let $V$ be an operator space, 
and $\mathbf{A}= (A_n)_{n=1}^\infty$ a well-founded matrix set over $V$.
Let $F = \spa\{ \delta_{\mathbf{A}}(x) \}_{x \in A_1} \subset \Fcb(\mathbf{A})$.
For any $m\in\N$ and $u \in M_m(F)$ we have
\[
\n{u}_{M_m(\Fcb(\mathbf{A}))} = \inf \bigg\{  \n{\sum_{\ell=1}^{L} \alpha_\ell(\alpha_\ell)^*}^{1/2} \n{\sum_{\ell=1}^{L} (\beta_\ell)^*\beta_\ell}^{1/2} \max_{1 \le \ell \le L} |c_\ell|\ \n{x^\ell - y^\ell}_{M_{n_\ell}} \bigg\}
\]
where the infimum is taken over all representations of $u$ of the form 
\begin{equation}\label{prop-Arens-Eells-incomplete-representations}
    u =  \sum_{\ell=1}^{L} \alpha_\ell c_\ell\big((\delta_{\mathbf{A}})_{n_\ell}(x^\ell) - (\delta_{\mathbf{A}})_{n_\ell}(y^\ell) \big) \beta_\ell
\end{equation}
where for each $1 \le \ell \le L$ 
we have $n_\ell \in\N$, $x^\ell,y^\ell \in A_{n_\ell}$, $\alpha_\ell \in M_{m,n_\ell}$, $\beta_\ell \in M_{n_\ell,m}$ ,$c_\ell \in \mathbb{K}$.
\end{proposition}

\begin{proof}
We first claim that any $u \in M_m(F)$ does indeed admit representations of the form given in the statement.
Note that for $x\in A_1$, the matrix in $M_m(F)$ having $\delta_{\mathbf{A}}(x)$ in the $(j,k)$ position and zeros elsewhere admits a representation of the desired form with $\alpha \cdot 1 (\delta_{\mathbf{A}}(x)-\delta_{\mathbf{A}}(0)) \beta$, with $\alpha = E_{j1} \in M_{m,1}$ and $\beta=E_{1k} \in M_{1,m}$ where $E_{rs}$ are the matrix units.
Since the matrices of the above form span all of $M_m(F)$, and it is clear that a representation for $u_1+u_2$ is obtained by adding up individual representations for $u_1$ and $u_2$, the claim is proved.

Now, use the infimum in the statement to define a seminorm $\tn{\cdot}_m$ on each space $M_m(F)$.
For a fixed $u \in M_m(F)$ with a representation as in \eqref{prop-Arens-Eells-incomplete-representations},
note that said representation is equivalent to having
\begin{equation}\label{prop-Arens-Eells-incomplete-representations-matrix-product}
    u = \begin{bmatrix}
\alpha_1 & \alpha_2 &\cdots &\alpha_L\\
\end{bmatrix} \cdot
\diag\Big( c_\ell\big( (\delta_{\mathbf{A}})_{n_\ell}x^\ell - (\delta_{\mathbf{A}})_{n_\ell}y^\ell\big)\Big)_{\ell=1}^L
\cdot \begin{bmatrix}
\beta_1 \\ \beta_2 \\ \vdots \\ \beta_L\\
\end{bmatrix}
\end{equation}
from where it is clear that $\n{u}_{M_m(\Fcb(\mathbf{A}))} \le \tn{u}_m$ by using Ruan's axioms for $\Fcb(\mathbf{A})$ together with Proposition \ref{prop-delta-is-cb-holomorphic}.  In particular, each $\tn{\cdot}_m$ is in fact a norm.
Additionally, note that for any $x,y\in A_m$ using a trivial representation we have $\tn{ (\delta_{\mathbf{A}})_m(x)-(\delta_{\mathbf{A}})_m(y) }_m \le \n{x-y}_{M_n(V)} = \n{ (\delta_{\mathbf{A}})_m(x)-(\delta_{\mathbf{A}})_m(y) }_{M_m(\Fcb(\mathbf{A}))}$, and therefore
\begin{equation}\label{tn-agrees-with-n-on-deltas}
    \tn{ (\delta_{\mathbf{A}})_m(x)-(\delta_{\mathbf{A}})_m(y) }_m  = \n{ (\delta_{\mathbf{A}})_m(x)-(\delta_{\mathbf{A}})_m(y) }_{M_m(\Fcb(\mathbf{A}))}.
\end{equation}

By using the same arguments as in \cite[Prop. 5.5]{chavez2026linearizing}, it is not difficult to verify that this sequence of norms $\tn{\cdot}_m$
satisfies the criteria in \cite[Prop. 2.3.6]{Effros-Ruan-book}, which then implies that $(\tn{\cdot}_n)_{n=1}^\infty$ defines a sequence of norms on $M_n(F)$ that satisfies Ruan's axioms,
and therefore defines an operator space structure on the completion of $F$ with respect to $\tn{\cdot}_1$. Let us denote this operator space by $\mathfrak{F}(\mathbf{A})$, and note that from the above calculations the formal identity on $F$ extends to a complete contraction $ \iota : \mathfrak{F}(\mathbf{A}) \to \Fcb(\mathbf{A})$,
and thus
 $\iota(\mathbf{B}_{\mathfrak{F}(\mathbf{A})}) \subseteq \mathbf{B}_{\Fcb(\mathbf{A})}$.
From \eqref{tn-agrees-with-n-on-deltas} and Proposition \ref{proposition-ball-is-convex-hull-of-image-of-delta}, we obtain $\overline{\iota(\mathbf{B}_{\mathfrak{F}(\mathbf{A})})}=\overline{\iota(\overline{\mathbf{B}}_{\mathfrak{F}(\mathbf{A})})} \supseteq \overline{\mathbf{B}}_{\Fcb(\mathbf{A})}$, and therefore $\overline{\iota(\mathbf{B}_{\mathfrak{F}(\mathbf{A})})} = \overline{\mathbf{B}}_{\Fcb(\mathbf{A})}$.
Standard arguments
(see e.g. \cite[Proof of Prop. 4.28]{Bowers-Kalton})
yield that $\iota(\overline{\mathbf{B}}_{\mathfrak{F}(\mathbf{A})}) \supseteq \mathbf{B}_{\Fcb(\mathbf{A})}$.
Multiplying by $\lambda\in(0,1)$ and taking the union, we get
$\iota(\mathbf{B}_{\mathfrak{F}(\mathbf{A})}) \supseteq \mathbf{B}_{\Fcb(\mathbf{A})}$
from where we obtain that 
$\iota(\mathbf{B}_{\mathfrak{F}(\mathbf{A})}) = \mathbf{B}_{\Fcb(\mathbf{A})}$
So, $\iota$ is not just a complete contraction but in fact a complete $1$-quotient, yielding the desired conclusion.
\end{proof}

We remark that in the representations \eqref{prop-Arens-Eells-incomplete-representations} we could omit the $c_\ell$'s, because they can be absorbed into the $\alpha_\ell$ or the $\beta_\ell$. However, including them is sometimes useful when doing calculations.
Note also that in the ``classical'' case where $\mathbf{A}$ is just $A_1$ (meaning that $A_n=\emptyset$ for $n\ge 2$), Proposition \ref{prop-Arens-Eells-incomplete-without-matrix} clearly reduces to \eqref{eqn-classical-AE-incomplete} in the case $m=1$, since $\alpha_\ell$ and $\beta_\ell$ become scalars.

Just as we related Lemma \ref{lemma-characterization-abs-mat-convex-hull} and Corollary \ref{cor-amconvex hull of increasing matrix form} by an equivalence between representations as sums and those as matrix products, the equivalence between representations as in \eqref{prop-Arens-Eells-incomplete-representations} or as in \eqref{prop-Arens-Eells-incomplete-representations-matrix-product}
immediately yields the following alternative statement for Proposition \ref{prop-Arens-Eells-incomplete-without-matrix}, which clearly generalizes
\cite[Thm. 5.2]{Braga-CD-Sinclair}.

\begin{corollary}\label{cor-Arens-Eells-incomplete}
Let $V$ be an operator space, 
and $\mathbf{A}= (A_n)_{n=1}^\infty$ a well-founded matrix set over $V$.
Let $F = \spa\{ \delta_{\mathbf{A}}(x) \}_{x \in A_1} \subset \Fcb(\mathbf{A})$.
For any $m\in\N$ and $u \in M_m(F)$ we have
\[
\n{u}_{M_m(\Fcb(\mathbf{A}))} = \inf \Big\{ \n{\alpha}_{M_{m,N}} \n{\beta}_{M_{N,m}} \max_{1 \le \ell \le L} |c_\ell| \n{x^\ell - y^\ell}_{M_{n_\ell}(V)} \Big\}
\]
where the infimum is taken over all representations of $u$ of the form $u = \alpha \cdot D \cdot \beta$ where $\alpha \in M_{m,N}$ and $\beta \in M_{N,m}$ are scalar matrices, and $D \in M_{N}(F)$  is a block diagonal matrix whose diagonal blocks are of the form $c_\ell((\delta_{\mathbf{A}})_{n_\ell}x^\ell-(\delta_{\mathbf{A}})_{n_\ell}y^\ell)$ with $c_\ell$ a scalar and $x^\ell,y^\ell \in A_{n_\ell}$ for  $1 \le \ell \le L$, with $n_1+\cdots+n_L=N$. 
\end{corollary}

For a general element of the classical Lipschitz-free space, its norm can be calculated through the nice formula \cite[Lemma 2.1]{AliagaPernecka2020}
\begin{equation}\label{eqn-norm-in-free-space-as-infinite-sum}
\n{u}_{\F(X)} = \inf \Big\{ \sum_{j=1}^\infty |c_j| d(x_j,y_j) \;:\; u = \sum_{j=1}^\infty c_j\big( \delta_X(x_j) - \delta_X(y_j) \big) \Big\}.
\end{equation}
This is a straightforward consequence of \eqref{eqn-classical-AE-incomplete}, together with the standard fact that if $Y$ is a dense subspace of a normed space $X$, then for every $x \in X$ and $\varepsilon>0$ we can find a norm convergent representation $x = \sum_{n=1}^\infty y_n$ with $y_n \in Y$ and $\sum_{n=1}^\infty \n{y_n} < \n{x}+\varepsilon$.
By the same general principle one can deduce the following result from Proposition \ref{prop-Arens-Eells-incomplete-without-matrix}, which we state without proof: one can simply translate the argument used to show Theorem \ref{thm-Arens-Eeels} below, via the aforementioned equivalence between representations as in \eqref{prop-Arens-Eells-incomplete-representations} or as in \eqref{prop-Arens-Eells-incomplete-representations-matrix-product}.

\begin{corollary}\label{cor-Arens-Eells-without-matrix}
Let $V$ be an operator space, 
and $\mathbf{A}= (A_n)_{n=1}^\infty$ a well-founded matrix set over $V$.
For any $m\in\N$ and $u \in M_m(\Fcb(\mathbf{A}))$ we have
\[
\n{u}_{M_m(\Fcb(\mathbf{A}))} = \inf \bigg\{  \n{\sum_{\ell=1}^{\infty} \alpha_\ell(\alpha_\ell)^*}^{1/2} \n{\sum_{\ell=1}^{\infty} (\beta_\ell)^*\beta_\ell}^{1/2} \sup_{\ell} |c_\ell|\ \n{x^\ell - y^\ell}_{M_{n_\ell}} \bigg\}
\]
where the infimum is taken over all representations of $u$ as a convergent series of the form 
\[
    u =  \sum_{\ell=1}^{\infty} \alpha_\ell c_\ell\big((\delta_{\mathbf{A}})_{n_\ell}(x^\ell) - (\delta_{\mathbf{A}})_{n_\ell}(y^\ell) \big) \beta_\ell
\]
where for each $\ell\in\N$ 
we have $n_\ell \in\N$, $x^\ell,y^\ell \in A_{n_\ell}$, $\alpha_\ell \in M_{m,n_\ell}$, $\beta_\ell \in M_{n_\ell,m}$, $c_\ell \in \mathbb{K}$.
\end{corollary}

To close out this section, we show how to use Corollary \ref{cor-Arens-Eells-incomplete} to obtain a description of the norms of arbitrary elements of $M_m(\Fcb(\mathbf{A}))$ in terms of products of infinite matrices.
Just as \eqref{eqn-norm-in-free-space-as-infinite-sum} is reminiscent of a similar way of calculating the norm of elements in the completion of the projective tensor product of Banach spaces \cite[Prop. 2.8]{MR1888309}, our next result should be compared to the description of the elements in the completion of the projective tensor product of operator spaces \cite[Thm. 10.2.1]{Effros-Ruan-book}.
We also point out that Theorem \ref{thm-Arens-Eeels} below generalizes \cite[Thm. 5.3]{Braga-CD-Sinclair}.

The notation for infinite matrices which we use is that of \cite{Effros-Ruan-book}. The reader can find all details in Sections 1.1, 10.1 and 10.2 of \cite{Effros-Ruan-book}, we will only recall some basics.
While \cite{Effros-Ruan-book} concentrates on the complex case, the real one is analogous.
Given an operator space $V$, $M_\infty(V)$ is the space of infinite matrices $(x_{jk})_{j,k=1}^\infty$ with entries $x_{jk}$ in $V$ whose truncations to $M_n(V)$ are uniformly bounded; $M_\infty(V)$ has a natural operator space structure.
A similar construction can be done for matrices $M_{I,J}(V)$, where $I$ (resp. $J$) is a finite or countably infinite set of indices for the rows (resp. columns); operator space structures and matrix multiplication can then be defined using truncations again (in this context, $\N$ is identified with $\infty$). Moreover, if $V=\K$, we simply write $M_\infty=M_\infty(\K)$, 
$M_{n,\infty}=M_{n,\infty}(\K)$, and $M_{\infty,n}=M_{\infty,n}(\K)$.

\begin{theorem}\label{thm-Arens-Eeels}
Let $V$ be an operator space, 
and $\mathbf{A}= (A_n)_{n=1}^\infty$ a well-founded matrix set over $V$.
For any $m\in\N$ and $u \in M_m(\Fcb(\mathbf{A}))$ we have
\[
\n{u}_{M_m(\Fcb(\mathbf{A}))} = \inf \Big\{ \n{\alpha} \n{\beta} \sup_{\ell\in\N} |c_\ell| \n{x^\ell - y^\ell} \Big\}
\]
where the infimum is taken over all representations of $u$ of the form $u = \alpha \cdot D \cdot \beta$ where $\alpha \in M_{m,\infty }$ and $\beta \in M_{\infty,m}$ are scalar matrices, and $D \in M_{\infty}(\spa\{\delta_{\mathbf{A}}(x)\}_{x \in A_1})$  is a block-diagonal matrix whose diagonal blocks are of the form $c_\ell\big( (\delta_{\mathbf{A}})_{n_\ell}x^\ell - (\delta_{\mathbf{A}})_{n_\ell}y^\ell \big)$ with $c_\ell$ a scalar and $x^\ell, y^\ell \in A_{n_\ell}$ for each $\ell \in \N$.
Furthermore, we can assume $\lim_{\ell\to\infty} |c_\ell| \n{x^\ell - y^\ell} = 0$.
\end{theorem}

\begin{proof}
 We start by proving the $\le$ inequality.
Suppose that we have a representation  $u = \alpha \cdot D \cdot \beta$ as in the statement.
By definition, this means that $u$ is the limit in $M_m(\Fcb(\mathbf{A}))$ of the truncated products $u(r) = \alpha^{r} \cdot D \cdot \beta^{r}$ as $r \to \infty$.
 By Corollary \ref{cor-Arens-Eells-incomplete} we get (by identifying the finitely supported matrix $\alpha^{r} \cdot D \cdot \beta^{r} = \alpha^{r} \cdot D^r \cdot \beta^{r}$ with a finite matrix) 
 if we take $N_r = n_1+n_2+\cdots n_r$,
\[
\n{u(N_r)}_{M_m(\Fcb(\mathbf{A}))} \le \n{\alpha^{N_r}} \n{\beta^{N_r}} \sup_{1 \le \ell \le r} |c_\ell| \n{x^\ell - y^\ell} 
\le \n{\alpha} \n{\beta} \sup_{\ell} |c_\ell| \n{x^\ell - y^\ell}  
\]
so taking the limit as $r\to\infty$ yields the $\le$ inequality.

Assume now without loss of generality that $u\not=0$ and $\n{u}_{M_m(\Fcb(\mathbf{A}))}<1$.
By density, there exists a sequence $\{u_s\}_{s=1}^\infty$ of nonzero terms in $M_m( \spa\{\delta_{\mathbf{A}}(x)\}_{x\in A_1} )$ such that
\[
u = \sum_{s=1}^\infty u_s \qquad \text{and} \qquad \sum_{s=1}^\infty \n{u_s}_{M_m(\Fcb(\mathbf{A}))} < 1.
\]
 Let $\varepsilon$ be a number with $0 < \varepsilon < 1 -\sum_{s=1}^\infty \n{u_s}_{M_m(\Fcb(\mathbf{A}))}$.
 For each $s$, Corollary \ref{cor-Arens-Eells-incomplete} gives  a representation $u_s = \alpha_s \cdot D_s \cdot \beta_s$, 
 with  $\alpha_s \in M_{m,N_s}$, $\beta_s \in M_{N_s,m}$, and $D_s \in M_{N_s}(\spa\{\delta_{\mathbf{A}}(x)\}_{x\in A_1} )$  a block diagonal matrix whose diagonal blocks are of the form $c_{\ell,s}\big((\delta_{\mathbf{A}})_{n_{\ell,s}}x^{\ell,s}-(\delta_{\mathbf{A}})_{n_{\ell,s}}y^{\ell,s}\big)$, $1 \le \ell \le L_s$,
 with $x^{\ell,s}, y^{\ell,s} \in A_{n_{\ell,s}}$
 and such that
 \[
 \n{\alpha_s} \n{\beta_s} \max_{1 \le \ell \le L_s} |c_{\ell,s}|\n{x^{\ell,s} - y^{\ell,s}}_{M_{n_{\ell,s}}(V)} < \n{u_s}_{M_m(\Fcb(\mathbf{A}))} + \varepsilon/2^s.
 \]
 Observe that for each $s \in \N$ we must have $\lambda_s=\max_{1 \le \ell \le L_s} |c_{\ell,s}|\n{x^{\ell,s} - y^{\ell,s}}_{M_{n_{\ell,s}}(V)}>0$ because $u_s \not=0$.
  Also, by rescaling we may assume
\[
\n{\alpha_s} = \n{\beta_s} < \bigg( \frac{ \n{u_s}_{M_m(\Fcb(\mathbf{A}))} + \varepsilon/2^s }{ \lambda_s } \bigg)^{1/2}
\]
Since
$
\sum_{s=1}^\infty \n{\alpha_s} \n{\beta_s} \lambda_s < 1,
$
we can find a sequence of numbers $p_s > 1$ converging to $\infty$ and such that 
$
\sum_{s=1}^\infty p_s \n{\alpha_s} \n{\beta_s} \lambda_s < 1.
$
 
 Define a block diagonal matrix $D \in M_\infty(\spa\{\delta_{\mathbf{A}}(x)\}_{x\in A_1})$ by concatenating the block diagonal matrices $p_s^{-1}\lambda_s^{-1} D_s$ along the diagonal of $D$.
 Note that, by design, $\n{D}<1$ and $\lim_{s\to\infty} \n{p_s^{-1}\lambda_s^{-1} D_s} = \lim_{s\to\infty} p_s^{-1} = 0$, so this matrix $D$ satisfies the condition in the ``furthermore'' part of the statement.
 Now define scalar matrices 
\[
\alpha =
\begin{bmatrix}
(p_1\lambda_1)^{1/2}\alpha_1 &(p_2\lambda_2)^{1/2}\alpha_2 &(p_3\lambda_3)^{1/2}\alpha_3 & \cdots\\
\end{bmatrix}
\in M_{m,\infty}
\]
 and
\[
\beta = \begin{bmatrix}
(p_1\lambda_1)^{1/2}\beta_1 &(p_2\lambda_2)^{1/2}\beta_2 & (p_3\lambda_3)^{1/2}\beta_3 &\cdots\\
\end{bmatrix}^T \in M_{\infty,m}.
\]
Clearly
\[
\n{\alpha}^2 = \n{\alpha \alpha^*} \le \sum_{s=1}^\infty p_s\lambda_s \n{\alpha_s}^2 = \sum_{s=1}^\infty p_s\lambda_s \n{\alpha_s} \n{\beta_s} <1, 
\]
so $\n{\alpha}<1$ and analogously $\n{\beta}<1$.
One can readily verify that $u = \alpha \cdot D \cdot \beta$, finishing the proof.
\end{proof}

\section{The linearization property}\label{sec-linearization-property}

As in the classical setting, completely Lipschitz functions linearize through their predual.

\begin{theorem}\label{thm-linearization}
Let $V$ be an operator space, 
and $\mathbf{A}= (A_n)_{n=1}^\infty$ a well-founded matrix set over $V$.
For every operator space $W$ and every completely Lipschitz function $f \in \Lipcb(\mathbf{A},W)$, there exists a unique linear map $\widehat{f}:\Fcb(\mathbf{A})\to W$ such that $\widehat{f}\circ \delta_{\mathbf{A}}=f$ and $\bn{\widehat{f}}_{\cb}=\n{f}_{ \Lipcb(\mathbf{A},W) }$. That is, $\widehat{f}$ makes the following diagram commutative:
\[
\xymatrix{\Fcb(\mathbf{A})\ar@{-->}[rd]^{\widehat{f}}\\
     A_1\ar[r]_{f}\ar^{\delta_{\mathbf{A}}}[u]& W }
\]
This property completely isometrically isomorphically characterizes $\Fcb(\mathbf{A})$.
Moreover, the assignment $f \mapsto \widehat{f}$ is a completely isometric isomorphism between $\Lipcb(\mathbf{A},W)$ and $\CB(\Fcb(\mathbf{A}),W)$.
\end{theorem}

\begin{proof}

We claim that one can define a linear map $ \widehat{f} : \spa \{ \delta_{\mathbf{A}}(x) \}_{x \in A_1} \to W$ by $\widehat{f}(\sum_{j=1}^n c_j \delta_{\mathbf{A}}(x_j)) = \sum_{j=1}^n c_j f(x_j)$ where $c_j \in \K$ and $x_j \in A_1$ for $j=1,\dotsc,n$.
If the evaluation functionals $\{ \delta_{\mathbf{A}}(x) \}_{x \in A_1 \setminus\{0\}}$ were linearly independent this would be obvious, but they may not be. Thus, we have to prove that two different representations of an element of $\spa \{ \delta_{\mathbf{A}}(x) \}_{x \in A_1}$ yield the same definition for $\widehat{f}$. Assume that $\sum_{j=1}^n c_j \delta_{\mathbf{A}}(x_j) = 0$ in $\Fcb(\mathbf{A})$. Therefore, for every $g \in \Lipcb(\mathbf{A}, \K)$ we have $0 = \pair{g}{\sum_{j=1}^n c_j \delta_{\mathbf{A}}(x_j)} = \sum_{j=1}^n c_jg(x_j)$.
In particular, for every $w' \in W'$ note that $w' \circ f \in \Lipcb(\mathbf{A},\K)$, and therefore
$0 = \sum_{j=1}^n c_j w'\circ f(x_j) = w' \Big( \sum_{j=1}^n c_j f(x_j) \Big)$, from where it follows that $\sum_{j=1}^n c_j f(x_j) =0$.
Therefore, $\widehat{f}$ is well defined on $\spa \{ \delta_{\mathbf{A}}(x) \}_{x \in A_1}$ and this is the only possible choice for $\widehat{f}$ on this subspace.

Now, fix $n\in\N$ and $u \in M_n(\spa \{ \delta_{\mathbf{A}}(x) \}_{x \in A_1})$.
For any $k\in\N$ and $w \in M_k(W')$, note that $\mpair{w'}{\widehat{f}_n(u)} = \mpair{w'\circ f}{u}$ where on the right we are identifying $M_k(W') \equiv \CB(W, M_k)$.
Therefore, by the duality between $\Fcb(\mathbf{A})$ and $\Lipcb(\mathbf{A})$,
\[
\n{\mpair{w'}{\widehat{f}_n(u)}}_{M_{nk}} \le \n{w'\circ f}_{M_k(\Lipcb(\mathbf{A}))} \n{u}_{M_n(\Fcb(\mathbf{A}))}
\le 
\n{w'}_{M_k(W')} \n{f}_{\Lipcb(\mathbf{A}),W)} \n{u}_{M_n(\Fcb(\mathbf{A}))}
\]
Taking supremum over all $w'$ of norm at most one then yields $\bn{\widehat{f}_n(u)}_{M_n(W)} \le \n{f}_{\Lipcb(\mathbf{A},W)} \n{u}_{M_n(\Fcb(\mathbf{A}))}$,
showing that $\widehat{f}$ is a completely bounded linear map on $\spa \{ \delta_{\mathbf{A}}(x) \}_{x \in A_1}$. Thus, $\widehat{f}$ uniquely extends to $\Fcb(\mathbf{A})$ with $\cb$-norm at most $\n{f}_{\Lipcb(\mathbf{A},W)}$.
The reverse inequality is obvious because $\widehat{f}$ extends $f$, so $\bn{\widehat{f}}_{\cb}=\n{f}_{ \Lipcb(\mathbf{A},W) }$.

The complete isometry conclusion now follows immediately from the identifications
\[
M_n(\CB(\Fcb(\mathbf{A}), W )) \equiv \CB(\Fcb(\mathbf{A}), M_n(W)), \qquad M_n(\Lipcb(\mathbf{A},W)) \equiv \Lipcb(\mathbf{A},M_n(W)).
\]

\end{proof}

As in the classical setting, we also have the following standard alternative linearization result:

\begin{theorem}\label{thm-cuadradito}
Let $V$ and $W$ be operator spaces, and $\mathbf{A}= (A_n)_{n=1}^\infty$, $\mathbf{B}= (B_n)_{n=1}^\infty$ well-founded matrix sets over $V$ and $W$, respectively.
Suppose that $f \in \Lipcb(\mathbf{A},W)$ is a function mapping $\mathbf{A}$ to $\mathbf{B}$, in the sense that for every $n\in \N$ and every $x \in A_n$ we have $f_n(x) \in B_n$.
Then there exists a unique 
completely bounded linear map $\widetilde{f} : \Fcb(\mathbf{A}) \to \Fcb(\mathbf{B})$ such that $\widetilde{f}\delta_{\mathbf{A}} = \delta_{\mathbf{B}}f$, that is, the following diagram commutes:
\[
\xymatrix{
\Fcb(\mathbf{A})\ar[r]^{\widetilde{f}} &\Fcb(\mathbf{B})\\
     A_1 \ar[r]^{f}\ar[u]^{\delta_{\mathbf{A}}}& B_1 \ar[u]^{\delta_{\mathbf{B}}}
     }
\]
Moreover, $\bn{\widetilde{f}}_{\cb} = \n{f}_{ \Lipcb(\mathbf{A},W)}$.

\end{theorem}

\begin{proof}
It follows immediately from Proposition \ref{prop-delta-is-cb-holomorphic} that
$\delta_{\mathbf{B}}f$ belongs to $\Lipcb(\mathbf{A},\Fcb(\mathbf{B}))$, 
and moreover $\n{\delta_{\mathbf{B}}f}_{\Lipcb(\mathbf{A},\Fcb(\mathbf{B}))} = \n{f}_{\Lipcb(\mathbf{A},W)}$. 
Thus, due to Theorem \ref{thm-linearization} the statement holds for $\widetilde f= \widehat{\delta_\mathbf{B}f}$.
Uniqueness is clear from the fact that the image of $\delta_{\mathbf{A}}$ has dense span in $\Fcb(\mathbf{A})$.

\end{proof}

\begin{remark}
For an operator space $V \not=\{0\}$, considering $B_V$ as a matrix set in the usual way (with all the levels above the first one being empty), we have $B_V \subseteq \mathbf{B}_V$. Taking the identity $B_V \to B_V$, Theorem \ref{thm-cuadradito} then provides a canonical  contraction $\F(B_{V}) \to \F^{\cb}(\mathbf B_{V})$.
Note that the latter map cannot be bijective, because in that case so would be its adjoint $\Lipcb(\mathbf{B}_V) \to \Lipo(B_V)$, which is just the inclusion map but it is not surjective by Examples \ref{examples completely Lip}.\eqref{examples completely Lip:modulo}.
However, for the simplest possible case $V=\R$, in the companion paper \cite{CD-Dimant-futuro} we show that $\F(B_{\R})$ and $\F^{\cb}(\mathbf B_{\R})$ are isomorphic as Banach spaces.
\end{remark}

\section*{AI disclosure statement}
Claude was used to study the converse of Corollary \ref{cor: N maximal complexification} (that is, to obtain Propositions \ref{prop-Vc-N-maximal-implies-V-2N-maximal} and \ref{prop-converse-N-maximal-fails}). The rest of the mathematical developments are stated, proved and verified by the authors, who take full responsibility for the content of the manuscript.
The authors also used AI tools for light language editing,  and reviewing for errata.

\bibliography{references}
\bibliographystyle{amsalpha}

\end{document}